\documentclass[aap,MSNbibl,citesort,dvips]{arximspdf}

\doi{10.1214/09-AAP622}
\volume{20}
\issue{2}
\pubyear{2010}
\firstpage{382}
\lastpage{429}

\makeatletter

\DeclareMathAlphabet\mathcaligr{OMS}{cmsy}{m}{n}

\newtheorem{theorem}{Theorem}[section]
\newtheorem{proposition}[theorem]{Proposition}
\newtheorem{lemma}[theorem]{Lemma}
\newtheorem{corollary}[theorem]{Corollary}

\newproclaim{Example}{Example}
\newproclaim{Remark}{Remark}

\newcommand{\xrightarrow}{\mathop{\hbox to 1cm{\rightarrowfill}}_}
\newcommand{\dd}{d}

\makeatother

\begin{document}
\begin{frontmatter}

\title{Asymptotic behavior of solutions of the fragmentation equation
with shattering: An approach via self-similar Markov processes}
\runtitle{Asymptotic behavior of solutions}

\begin{aug}
\author[A]{\fnms{B\'en\'edicte} \snm{Haas}\corref{}\ead[label=e1]{haas@ceremade.dauphine.fr}}
\runauthor{B. Haas}
\affiliation{Universit\'{e} Paris-Dauphine}
\address[A]{Universit\'{e} Paris-Dauphine\\
Place du Mar\'{e}chal\\
\quad de Lattre de Tassigny\\
Paris, 75775\\
France\\
\printead{e1}} 
\end{aug}

\pdfauthor{Benedicte Haas}

\received{\smonth{12} \syear{2008}}
\revised{\smonth{6} \syear{2009}}

%
\begin{abstract}
The subject of this paper is a fragmentation equation with
nonconservative solutions, some mass being lost to a dust of zero-mass
particles as a consequence of an intensive splitting. Under some
assumptions of regular variation on the fragmentation rate, we describe
the large time behavior of solutions. Our approach is based on
probabilistic tools: the solutions to the fragmentation equation are
constructed via nonincreasing self-similar Markov processes that
continuously reach 0 in finite time. Our main probabilistic result
describes the asymptotic behavior of these processes conditioned on
nonextinction and is then used for the solutions to the fragmentation
equation.

We note that two parameters significantly influence these large time
behaviors: the rate of formation of ``nearly-1 relative masses'' (this
rate is related to the behavior near $0$ of the L\'evy measure
associated with the corresponding self-similar Markov process) and the
distribution of large initial particles. Correctly rescaled, the
solutions then converge to a nontrivial limit which is related to the
quasi-stationary solutions of the equation. Besides, these
quasi-stationary solutions, or, equivalently, the quasi-stationary
distributions of the self-similar Markov processes, are fully
described.
\end{abstract}

%

\end{frontmatter}

\section{Introduction and main results}\label{sec1}

Fragmentation processes occur in a variety of natural phenomena,
including polymer degradation, mineral grinding and droplet break-up,
but also in the analysis of algorithms, phylogeny, etc.
The kinetic equation used in the physics literature to describe the
time-evolution of masses of particles prone to fragmentation has the form
%
%
\begin{equation}
\label{usualequation}
\partial_t n_t(x)=\int_x^{\infty} a(y)b(y,x)n_t(y) \,dy -a(x)n_t(x),
\end{equation}
where $n_t(x)$ is the concentration of particles of mass $x$ at time
$t$, $a(x)$ is the overall rate at which a particle with mass $x$
splits and $b(x,y)$ describes the distribution of particles of mass $y$
produced by the fragmentation of a particle of mass $x$. It is assumed
that no mass is lost when a particle breaks up, that is, $\int_0^x y
b(x,y) \, dy =x$.
The integral in the right-hand side of (\ref{usualequation}) models
the increase of particles of mass $x$ due to the fragmentation of
particles of masses $y>x$, whereas the negative term $-a(x)n_t(x)$
models the loss of particles of mass $x$, due to their fragmentation
into smaller particles. This fragmentation equation has been
intensively studied by both physicists and mathematicians. Among the
first papers on the topic, we may cite, for example, \cite{McGZ,Melzak}.

In both the physics and mathematics literature, particular attention
has been paid to models with the following self-similar dynamic:
\begin{itemize}
\item $a(x)=Cx^{\alpha}$, for some fixed $C>0$ and $\alpha\in
\mathbb R$;
\item $b(x,y)=h(y/x)/x$ [with $h$ such that $\int_0^1 uh(u)\,
du=1$]. This means that the distribution of the ratios of daughter
masses to parent mass is only determined by a function of these ratios
(and not by the parent mass).
\end{itemize}
There are two reasons for this. These self-similar assumptions are
relevant for applications, for example, for polymer degradation
\cite{ZMcG}, mineral crushing in the mining industry (\cite{BMEnergy} and
the references therein) and the construction of phylogenetic trees
\cite{Aldous96}. But they are also more mathematically tractable.
For the same reasons, there is also a significant literature on
probabilistic models for the microscopic mechanism of fragmentation
with a self-similar dynamic. We refer to the book by Bertoin \cite
{BertoinBook} for an overview and to the papers \cite{ESRR} and \cite
{HaasLossMass} for discussions of the relations between the
probabilistic models and the above equation.

The goal of this paper is to contribute to the understanding of
solutions of the self-similar fragmentation equation, by describing
their large time behavior. The cases where $\alpha>0$ are treated in
\cite{ESRR} and we will be concerned here only with the negative cases
$\alpha<0$.

We will actually consider the following generalization of the weak form
of the above fragmentation equation (\ref{usualequation}) with a
self-similar dynamic:
%
%
\begin{equation}
\label{eqfrag}
\partial_t\langle\mu_t,f\rangle
=\int_{0}^{\infty} x^{\alpha} \biggl( \int_{0}^1
\bigl(f(yx)-f(x)y \bigr) B(dy) \biggr) \mu_t( dx),
\end{equation}
where $(\mu_t,t\geq0)$ denotes a family of measures on $]0,\infty[$,
$\alpha\in\mathbb R$, $B$ is a measure on $]0,1[$ such that
%
%
\begin{equation}
\label{contraintesB}
\int_0^{1}y (1-y)B( dy)<\infty\quad\mbox{and}\quad B (]0,1[)>0,
\end{equation}
and $f$ denotes any test function. When $B( dy)=Ch(y) \,
dy$ with $\int_0^1 yh(y)\, dy=1$ and $\mu_t( dx) =
n_t(x)\, dx$, we recover the weak form of (\ref{usualequation})
with $a(x)=Cx^{\alpha}$ and $b(x,y)=h(y/x)/x$. Informally, (\ref
{eqfrag}) corresponds to models in which particles with mass $xy$,
$0<y<1$, are produced from the splitting of a particle with mass $x$ at
rate $x^{\alpha}B( dy)$. Note that the overall rate at which a
particle with mass $x$ splits is $x^{\alpha} \int_0^1 y B( d
y)$, which may be infinite here. Let us add that the physical
interpretation of the fragmentation equation imposes some constraints
on the measure $B$. However, other interpretations are possible and, in
the following, we will be concerned with all measures $B$ satisfying
(\ref{contraintesB}).

We focus on solutions of (\ref{eqfrag}) with finite and nonzero
initial total mass. The fragmentation equation being linear, we
suppose, without loss of generality, that $\int_{0}^{\infty} x \mu_0
( d x)=1$. To be precise, we call \textit{a solution of} (\ref
{eqfrag}) \textit{starting from $\mu_0$} any family of measures $(\mu_t,t\geq
0)$ on $]0,\infty[$ starting from $\mu_0$ such that:

\begin{itemize}
\item $(\mu_t,t\geq0)$ satisfies (\ref{eqfrag}) for any
test function $f \in C^1_c$, the set of real-valued continuously
differentiable functions on $]0,\infty[$ with compact support;
\item the natural ``physical properties''
%
%
\begin{equation}
\label{gainmass}
m(t):=\langle\mu_t,\mathrm{id}\rangle\leq m(0)=1\qquad \forall t \geq0,
\end{equation}
and
%
%
\begin{equation}
\label{creationmass}
\mu_0( [M,\infty[)=0 \qquad\mbox{for some }M>0 \quad\Rightarrow\quad\mu_t
([M,\infty[)=0 \qquad \forall t \geq0,\hspace*{-12pt}
\end{equation}
are respected (``$\mathrm{id}$'' denotes the identity function).
\end{itemize}

Note the \textit{self-similarity of solutions}: if $(\mu_t,t \geq0)$
is a solution of (\ref{eqfrag}), then so is $(\gamma^{-1} \mu
_{t\gamma^{\alpha}}\circ(\gamma\,\mathrm{id})^{-1})$ for all $\gamma
>0$. Also, note that if $(\mu_t,t \geq0)$ is a solution of the
equation with parameters $(\alpha,B)$, then for all $c>0$, $(\mu
_{ct},t \geq0)$ is a solution of the equation (\ref{eqfrag}) with
parameters $(\alpha,cB)$.

Many results on the existence and uniqueness of solutions of (\ref
{usualequation}) are available in the literature; see, for example,
\cite{Banasiaksurvey,ESRR,LaurencotCoagFrag} and the references
therein. With the definition above, we have the following result on the
existence and uniqueness of solutions of (\ref{eqfrag}), which is a
generalization of Theorem 1 of \cite{HaasLossMass} (see also \cite
{FG} for a similar approach). We recall that a \textit{subordinator} is
a nondecreasing L\'evy process and that its distribution is
characterized by two parameters: a nonnegative drift coefficient and a
so-called L\'evy measure on $]0,\infty[$ that governs the jumps of the
process. See Section \ref{background} for background on this topic.
\begin{theorem}
\label{theu}
Let $\mu_0$ be a measure on $]0,\infty[$ such that $\int_{0}^{\infty
} x \mu_0( d x)=1$ and let $\xi$ be a subordinator with zero
drift and L\'evy measure $\Pi$ given, for any measurable function
$g\dvtx
]0,\infty[\  \rightarrow[0,\infty[$, by
%
%
\begin{equation}
\label{PiB}
\int_0^{\infty}g(x)\Pi( d x)=\int_0^{1}g(-\ln(x))x
B( d x).
\end{equation}
Then, for each $t >0$, define a measure $\mu_t$ on $]0,\infty[$ by
%
%
\begin{equation}
\label{solution}
\int_0^{\infty}f(x) x\mu_t( dx):=\int_0 ^{\infty} \mathbb
E \bigl[f \bigl(x\exp\bigl(-\xi_{\rho(x^{\alpha}t)}\bigr) \bigr) \bigr]
x\mu_0( d x)
\end{equation}
for all measurable $f \dvtx[0,\infty[ \ \rightarrow[0,\infty[$, $f(0)=0$,
where $\rho$ is the time-change
\[
\rho(s):= \inf\biggl\{u \geq0 \dvtx\int_0^u\exp(\alpha\xi_r)
\, d r >s \biggr\}\qquad \forall s \geq0.
\]

\textup{(i)} The family $(\mu_t,t\geq0)$ is a solution of (\ref
{eqfrag}), provided that
\[
\alpha\leq0
\]
or
\[
\alpha>0 \quad\mbox{and either}\quad \int_{1}^{\infty} x\ln(x) \mu
_0( d x)<\infty
\]
(or)
\[
x \in\,]0,1[ \ \rightarrow x^{|\alpha
|} \int_0^{x}y B( dy) \mbox{ is bounded near 0}.
\]

\textup{(ii)} This solution is unique, provided that $\mu_0([M,\infty
[)=0$ for some $M>0$.
\end{theorem}

When the family $(\mu_t,t \geq0)$ is constructed via a subordinator
by (\ref{solution}), some conditions on $\mu_0$ and $B$ for the
existence of a density for $\mu_t$, $t > 0$, can be stated explicitly;
see, for example, \cite{FGdensities}, Proposition 3.10. We also recall
that there may exist multiple solutions of the fragmentation equation
when the assumption (\ref{gainmass}) is dropped. We refer to \cite
{Banasiaksurvey} for some explicit examples.

The proof of Theorem \ref{theu}, based on that of Theorem 1 in \cite
{HaasLossMass}, is postponed to the \hyperref[App]{Appendix}.

The main purpose of this paper is to use the construction (\ref
{solution}) of solutions of the fragmentation equation to describe the
large time behavior of these solutions when $\alpha<0$. From another,
but equivalent, point of view, our main results describe the large time
behavior of exponentials of minus time-changed subordinators, as
defined in Theorem \ref{theu}, conditioned on nonextinction. These
processes belong to the family of so-called \textit{self-similar
Markov processes}. We refer to Section \ref{SectionSub} for a
statement of our results in that context.

The study of the large time behavior of solutions of the fragmentation
equation when $\alpha> 0$ is investigated in detail in \cite{ESRR}.
We point out that some results of \cite{ESRR} can be redemonstrated
using a probabilistic approach: it consists mainly of combining the
subordinator construction of solutions of the fragmentation equation
with the description of large time behavior of time-changed
subordinators when $\alpha> 0$ investigated in \cite{BC}.

From now on, we consider $\alpha<0$. It is well known that in such a
case, small particles split so quickly that they are reduced to a dust
of zero-mass particles, so that the total mass of nonzero particles
\[
m(t)=\langle\mu_t, \mathrm{id}\rangle
\]
decreases as time passes. This phenomenon, sometimes called
``shattering,'' has been studied in, for example, \cite
{Banasiaksurvey,BertoinAB,Filippov,HaasLossMass,McGZ,Wagner}. More
precisely, one can check that the total mass $m$ is strictly decreasing
and strictly positive on $[0,\infty[$, and that $m(t) \rightarrow0$
as $t \rightarrow\infty$; see the forthcoming Proposition \ref
{massgeneral} for a proof in our framework.

In order to describe the behavior of $m(t)$ as $t \rightarrow\infty$
more accurately,
we introduce the function defined for all $t \geq0$ by
%
%
\begin{equation}
\label{defphi}
\phi(t):=\int_0^1 (1 - x^{t} )xB( d x).
\end{equation}
It is not hard to check that the function $t \rightarrow t/\phi(t)$ is
continuous and strictly increasing on $]0,\infty[$, and that its range
is $](\int_0^1 |{\ln(x)}| x B( d x) )^{-1},\infty[$. Note that
the integral $\int_0^1 |{\ln(x)}| x B( d x)$ may be finite or
infinite. Then, introduce
%
%
\begin{equation}
\varphi,\qquad \mbox{the inverse of }
t \rightarrow t/\phi(t),
\end{equation}
which is well defined in a neighborhood of $\infty$. This function
will play a key role in the description of the long-time behavior of
solutions of the fragmentation equation.

Most of our main results rely on the following hypothesis on the
measure $B$:
\renewcommand{\theequation}{H}
\begin{eqnarray}
\label{hypothese}
&&\mbox{the function }u \dvtx]0,1[\  \rightarrow\int_0^{1-u}x B( d
x) \nonumber\\[-12pt]\\[-12pt]
&&\mbox{varies regularly at } 0 \mbox{ with an index }
{-}\beta\in\,]{-1},0],\nonumber
\end{eqnarray}
which, in particular, ensures that $\phi$ and $\varphi$ are regularly
varying functions at $\infty$ with respective indices $\beta$ and
$1/(1-\beta)$. See Section \ref{secRV} for details and background on
regular variation.

Finally, we mention that the large time behavior of solutions of the
fragmentation equation will depend strongly on the structure of the
initial measure $\mu_0$, mainly on the manner in which it distributes
weight near $\infty$. The statements of our results are therefore
split into two parts, according as to whether the initial measure has a
bounded support (Section \ref{compact}) or not (Section \ref
{noncompact}). Section \ref{sssolutions} deals with the
quasi-stationary solutions.

\subsection{Initial measure $\mu_0$ with bounded support}
\label{compact}

In this subsection, we adopt the following hypotheses and notation:

\begin{itemize}
\item $\alpha<0$;
\item the measure $\mu_0$ has a bounded support, that is, $\mu
_0([M,\infty[)=0$ for some \mbox{$M>0$};
\item $(\mu_t,t \geq0)$ denotes the unique solution of the fragmentation
equation (\ref{eqfrag}) starting from $\mu_0$.
\end{itemize}

The supremum of the support of $\mu_0$ is the real number
$s$ such that $\mu_0(]s,\infty[)=0$ and $\mu_0(]s-\varepsilon,s
])>0$ for all $\varepsilon<s$.
Thanks to the self-similarity of solutions, we can, and will, always
suppose that \textit{this supremum is equal to $1$}.
In such a framework, we have the following results.
\begin{proposition}
\label{mproperties}
For all $\lambda<\phi(\infty):=\lim_{x \rightarrow\infty} \phi
(x)$, there exists a constant $C_{\lambda}<\infty$ such that
\[
m(t) \leq C_{\lambda}\exp(-\lambda t)\qquad \forall t \geq0.
\]
More precisely,
under the hypothesis (\ref{hypothese}),
\[
-\ln(m(t)) \mathop{\sim}_{t \rightarrow\infty} \frac{(1-\beta
)}{|\alpha|}
\varphi(|\alpha|t).
\]
In particular, $t \rightarrow-\ln(m(t))$ is regularly varying at
$\infty$ with index $1/(1-\beta)$.
\end{proposition}

Together with the following theorem, this gives a complete description
of the large time behavior of $(\mu_t, t \geq0)$. Here, two positive
functions $g$ and $h$ are said to be asymptotically equivalent if
$g(x)/h(x)\rightarrow1$ as $x \rightarrow\infty$.
\begin{theorem}
\label{theq}
Suppose that (\ref{hypothese}) holds and $\int_0 |{\ln(x)}|xB(
dx)<\infty$. Then, for all continuous bounded test functions $f\dvtx
]0,\infty[\  \rightarrow\mathbb R$,
\[
\frac{1}{m(t)} \int_0^{\infty} f \biggl( \biggl( \frac{\varphi
(|\alpha|t)}{|\alpha|t} \biggr)^{1/|\alpha|}x
\biggr)x\mu_t( d x) \mathop{\rightarrow}_{t \rightarrow
\infty} \int_0^{\infty} f(x) x\mu_{\infty}( d x),
\]
where $x\mu_{\infty}( d x)$ is a probability distribution on
$]0,\infty[$ that is characterized by its moments
%
%
\setcounter{equation}{9}
\renewcommand{\theequation}{\arabic{equation}}
\begin{equation}
\label{muinfinie}
\int_0^{\infty} x^{|\alpha| n}x\mu_{\infty}( d x)=\phi
(|\alpha|)\phi(2|\alpha|)\cdots\phi(n|\alpha|),\qquad n \geq1.
\end{equation}
The function $t \rightarrow\varphi(|\alpha|t)/(|\alpha|t)$ can be
replaced by any asymptotically equivalent function.
\end{theorem}

It is interesting to compare this result with that obtained by
Escobedo, Mischler and Rodriguez Ricard
\cite{ESRR} when the parameter $\alpha$ is positive. As
already mentioned, part of their result can be rediscovered and
completed by using results of Bertoin and Caballero \cite{BC}. With
our notation, and under the assumptions $\int_0 |{\ln(x)}| x B(
dx)<\infty$ and $\alpha>0$, the asymptotic behavior of the solution
$(\mu_t,t\geq0)$ of the fragmentation equation $(\alpha,B)$ starting
from $\mu_0=\delta_1$ can be described as follows:
\[
\int_0^{\infty} f (t^{1/\alpha}x
)x\mu_t( d x) \mathop{\rightarrow}_{t \rightarrow
\infty} \int_0^{\infty} f(x) x\eta_{\infty}( d x)
\]
for all continuous bounded functions $f\dvtx]0,\infty[\  \rightarrow
\mathbb R$. The measure $x\eta_{\infty}( d x)$ is a
probability measure on $]0,\infty[$. Interestingly, the measure $B$ is
then involved only in the description of the limit measure $\eta
_{\infty}$, not in the ``shape'' of the speed of decrease of masses to $0$.

We return to the case $\alpha<0$. Note that when $\int
_0^{1-u}xB( dx) \sim u^{-\beta}$ as $u\rightarrow0$ for some
$\beta\in[0,1[$, we have $\phi(t) \sim\Gamma(1-\beta)t^{\beta}$
and therefore $(\varphi(|\alpha|t)/|\alpha|t )^{1/|\alpha|} \sim
C_{\alpha,\beta} t^{\beta/((1-\beta)|\alpha|)}$ as $t \rightarrow
\infty$, where $C_{\alpha,\beta}=(|\alpha|^{\beta}\Gamma(1-\beta
))^{1/((1-\beta)|\alpha|)}$. When we also have $\int_0 |{\ln(x)}|xB( dx)<\infty$,
Theorem \ref{theq} then reads
\[
\frac{1}{m(t)} \int_0^{\infty} f \bigl( C_{\alpha,\beta} t^{\beta
/((1-\beta)|\alpha|)} x
\bigr)x\mu_t( d x) \mathop{\rightarrow}_{t \rightarrow
\infty} \int_0^{\infty} f(x) x\mu_{\infty}( d x)
\]
for all continuous bounded test functions $f\dvtx]0,\infty[\  \rightarrow
\mathbb R$.

The existence and uniqueness of a measure $\mu_{\infty}$ on
$]0,\infty[$ satisfying (\ref{muinfinie}) actually hold without any
assumption of regular variation on the measure $B$ or assumptions on
its behavior near $0$; see the discussion near equation (\ref{defmuR})
in Section \ref{SectionSub} for details. Some properties of the
measure $\mu_{\infty}$ (tail behavior near $0$ and near $\infty$)
are given in Section \ref{PropertiesZ}. In Section \ref{sssolutions},
we discuss its links with
the \textit{quasi-stationary solutions} of the fragmentation equation.

The proof of Theorem \ref{theq} consists of describing the behavior of
the mass of a typical random nondust particle, defined as follows: at
each time $t$, choose a particle at random among the particles with a
\textit{strictly positive mass}, with a probability proportional to
its mass. That is, if $M(t)$ denotes the mass of this random particle,
then the distribution of $M(t)$ is given by
\[
M(t) \stackrel{\dd}{\sim} \frac{x\mu_t( d x)}{m(t)}.
\]
In other words, in terms of the subordinator $\xi$ related to the
equation by (\ref{solution}),
$M(t)$ is distributed as $M(0)\exp(-\xi_{\rho(M(0)^{\alpha}t)})$,
conditioned to be strictly positive, with $M(0)$ independent of $\xi$.
In terms of $M$, the statement of Theorem \ref{theq} can be rephrased
as follows:
\[
\biggl( \frac{\varphi(|\alpha|t)}{|\alpha|t} \biggr) ^{1/|\alpha
|} M(t) \stackrel{\dd}{\rightarrow}
M_{\infty},
\]
where $ M_{\infty}$ is a random variable with distribution $x\mu
_{\infty}( dx)$. Note the special case $\int_0^{1}x B(
d x)<\infty$, where $\varphi(t)/t \rightarrow\int_0^{1}x B(
d x)<\infty$. We then have that $M(t)$ converges in distribution to a
nontrivial limit. In the other\vspace*{-2pt} cases satisfying the assumptions of
Theorem \ref{theq}, $\varphi(t)/t \rightarrow\infty$ and therefore
$M(t) \stackrel{\mathbb{P}}{\rightarrow}0$.

Using this random approach, we can also specify the behavior of masses
that decrease at different speeds to $0$, as follows.
\begin{proposition}
\label{fasterslower}
Assume that (\ref{hypothese}) holds and let $\kappa:=\int_0^1 |{\ln(x)}|xB( dx)<\infty$.

\begin{longlist}
\item Suppose, moreover, that the support of $B$ is not
included in a set of the form $\{ a^n, n \in\mathbb N\}$ for some $a
\in\,]0,1[$. Then, for all measurable functions $g \dvtx[0,\infty[
\ \rightarrow\,]0,\infty[$ converging to $0$ at $\infty$,
\[
\frac{g(t)^{\alpha}}{m(t)} \int_0^{g(t) (\varphi(|\alpha
|t)/|\alpha|t)^{1/\alpha}} x \mu_t ( dx) \mathop{\rightarrow
}_{t \rightarrow\infty} \frac{1}{|\alpha|\kappa}.
\]
\item For all measurable functions $g \dvtx[0,\infty[
\ \rightarrow\,]0,\infty[$ converging to $\infty$ at $\infty$:
\begin{itemize}
\item if $g^{|\alpha|}(t)t/\varphi(t)$ converges to $\infty$ at
$\infty$,
then
\[
\int_{g(t) (\varphi(|\alpha|t)/|\alpha|t)^{1/\alpha}} ^{\infty} x
\mu_t ( dx)=0
\]
for all $t$ sufficiently large;
\item if $g^{|\alpha|}(t)t/\varphi(t)$ converges to $0$ at
$\infty$ and $0<\beta<1$, then
\[
\limsup_{t \rightarrow\infty} \frac{1} {\phi^{-1}(g(t)^{|\alpha
|}) }\ln\biggl( \frac{\int_{g(t) (\varphi(|\alpha|t)/|\alpha
|t)^{1/\alpha}} ^{\infty} x \mu_t ( dx)}{m(t)} \biggr) \leq
-\frac{\beta}{|\alpha|},
\]
where $\phi^{-1}$ denotes the inverse of $\phi$.
\end{itemize}
\end{longlist}
\end{proposition}

Note that the first assertion of (ii) is obvious since $g(t)
(\varphi(|\alpha|t)/|\alpha|t)^{1/\alpha} \rightarrow\infty$
(which means that for $t$ sufficiently large, it is larger than $1$,
the supremum of the support of $\mu_t$).

We conclude this section with the following result on the remaining
mass at time $t$ of particles of mass $1$ when $\mu_0(\{ 1\})>0$. The
measure $\mu_{\infty}$ is that introduced in Theorem \ref{theq}.
\begin{proposition}
\label{mass1}
Suppose that $\mu_0(\{ 1\})>0$ and set $\phi(\infty):=\int_0^1 x
B( d x) \in\,]0,\infty]$. Then, for all $t \geq0$,
\[
\mu_t(\{ 1\})=\exp(-t\phi(\infty) )\mu_0(\{ 1\}).
\]
If, further, (\ref{hypothese}) is satisfied, $\int_0 |{\ln(x)}| x
B( dx)<\infty$ and $\phi(\infty)<\infty$, then
\[
\frac{\mu_t(\{ 1\})}{m(t)} \mathop{\rightarrow}_{t \rightarrow
\infty} \phi(\infty)^{1/|\alpha|} \mu_{\infty} \bigl( \bigl\{\phi
(\infty)^{1/|\alpha|} \bigr\} \bigr)
\]
and this limit is nonzero if and only if
$\int^{1} \frac{B( d x)}{1-x}<\infty.
$
\end{proposition}

This means that under the assumptions of Proposition \ref{mass1}, for
large times, the remaining total mass of mass-$1$ particles is
proportional to the total mass of nonzero particles when $\int^{1}
(1-x)^{-1} B( d x)<\infty$, whereas it is negligible compared
to the total mass of nonzero particles when $=\int^{1} (1-x)^{-1}
B( d x)=\infty$.
We point out that the convergence of Proposition \ref{mass1} is
\textit{not} necessarily true when $\mu_0(\{ 1\})=0$ [since then $\mu
_t(\{ 1\})=0$ for all $t \geq0$, whereas the term in the limit may be
strictly positive].

\subsection{Initial measure $\mu_0$ with unbounded support}
\label{noncompact}

We still suppose that \mbox{$\alpha<0$} and we denote by $(\mu_t,t \geq0)$
the solution of the fragmentation equation (\ref{eqfrag}) starting
from $\mu_0$ and constructed via a subordinator by formula (\ref{solution}).
The asymptotic behavior of the mass $m(t)$ is then strongly modified by
the presence of large masses and depends on the behavior as $t
\rightarrow\infty$ of both $\phi(t)$ and $ \mu_0 ( [t,\infty
[ )$. We investigate two particular cases: exponential and power
decreases of $ \mu_0 ( [t,\infty[ )$ as $t \rightarrow
\infty$.
\begin{theorem}
\label{notcompact} Assume that (\ref{hypothese}) holds and that $\mu
_0$ possesses a density, say~$u_0$, in a neighborhood of $\infty$ such that
\[
\ln( u_0(x) ) \mathop{\sim}_{\infty}-Cx^{\gamma}
\]
for some $\gamma>0$.

\begin{longlist}
\item Then,
\[
-\ln( m(t) ) \mathop{\sim}_{\infty} C_{\alpha,\beta
,\gamma} C^{(1+(1-\beta)\gamma/|\alpha|)^{-1}} h(t),
\]
where $h$ is the inverse, well defined in the neighborhood of $\infty
$, of $t \rightarrow t^{1+|\alpha|/\gamma}/\phi(t)$ and
\[
C_{\alpha,\beta,\gamma}= \bigl(1+|\alpha|^{-1}\gamma(1-\beta)
\bigr) \biggl( \frac{|\alpha|^{1/(1-\beta)}}{\gamma}
\biggr)^{{\gamma(1-\beta)}/({\gamma(1-\beta)+|\alpha|})}.
\]
In particular, $-\ln(m(t))$ varies regularly at $\infty$ with index
$1/(1-\beta+|\alpha|/\gamma)$.

\item Suppose, moreover, that $\int_0 |{\ln(x)}|xB(
dx)<\infty$, which ensures that the function $\ln(m)$ is
differentiable on $]0,\infty[$. Then, if the derivative $(\ln(m))'$
is regularly varying at $\infty$, one has, for all continuous bounded
test functions $f\dvtx]0,\infty[\  \rightarrow\mathbb R$,
\[
\frac{1}{m(t)} \int_0^{\infty} f \biggl( \biggl( \frac
{h(t)}{C_{\alpha,\beta,\gamma,C}t} \biggr)^{1/|\alpha|}x
\biggr)x\mu_t( d x) \mathop{\rightarrow}_{t \rightarrow
\infty} \int_0^{\infty} f(x) x\mu_{\infty}( d x),
\]
where $\mu_{\infty}( d x)$ is the measure introduced in
Theorem \ref{theq} and
\[
C_{\alpha,\beta,\gamma,C}=\frac{C_{\alpha,\beta,\gamma}
C^{(1+(1-\beta)\gamma/|\alpha|)^{-1}}}{1-\beta+|\alpha|/\gamma}.
\]
\end{longlist}
\end{theorem}

Assuming that the derivative $(\ln(m))'$ is regularly varying at
$\infty$ may seem overly demanding. In actual fact, this assumption is
also needed to obtain Theorem \ref{theq}, but we are able to show that
it is always satisfied under the hypotheses of this theorem (see Lemma
\ref{tech3}). Unfortunately, it seems difficult to adapt this proof to
the case where the measure $\mu_0$ has unbounded support. However,
according to a classical result on regular variation (the monotone
density theorem), $(\ln(m))'$ varies regularly at $\infty$ provided
that $\ln(m)$ varies regularly at $\infty$ and $(\ln(m))'$ is
monotone near $\infty$, which can be checked in some particular cases.

There is also the following result on the decrease of the mass $m$ when
the density $u_0$ of $\mu_0$ has a power decrease near $\infty$.
\begin{proposition}
\label{powernc}
Assume that $\mu_0$ possesses a density $u_0$ in a neighborhood of
$\infty$ such that
\[
u_0 (x) \mathop{\sim}_{\infty}Cx^{-\gamma}
\]
for some $\gamma>2$. Then,
\[
m(t) \mathop{\sim}_{\infty} C' t^{({\gamma-2})/{\alpha}}
\]
with $C'=|\alpha|^{-1}C \int_0^{\infty}\overline m(u)u^{
({2-\gamma})/{\alpha}-1} \, du<\infty$, where $\overline m$
denotes the total mass of the solution of the fragmentation equation
with the same parameters $\alpha, B$ as that considered here, and with
initial distribution $\delta_1$, the Dirac mass at $1$.
\end{proposition}

\subsection{Quasi-stationary solutions}
\label{sssolutions}
A \textit{quasi-stationary solution} of the fragmentation equation
(\ref{eqfrag}) is a solution $(\mu_t,t \geq0)$ such that
\[
\mu_t=m(t)\mu_0\qquad \forall t \geq0,
\]
with $m(t)=\langle\mu_t, \mathrm{id}\rangle$. These quasi-stationary solutions
are closely related to the measure $\mu_{\infty}$ introduced in the
statement of Theorem \ref{theq}. We have already mentioned that
existence and uniqueness of such a measure $\mu_{\infty}$ satisfying
(\ref{muinfinie}) hold without any assumption of regular variation on
the measure $B$ or on its behavior near $0$. The interesting fact is
that, whatever the conditions on $B$, this measure and its self-similar
counterparts
\[
\mu^{(\lambda)}_{\infty}:=\lambda^{-1}\mu_{\infty}\circ(\lambda\,
\mathrm{id})^{-1},
\]
$\lambda>0$, are the only initial measures leading to quasi-stationary
solutions of the fragmentation equation (\ref{eqfrag}).
\begin{theorem}
\label{stat}
For all $\lambda>0$, let $(\mu^{(\lambda)}_{\infty,t}, t \geq0 )$
denote the solution of the fragmentation equation (\ref{eqfrag})
starting from $\mu^{(\lambda)}_{\infty}$ and constructed via a
subordinator by (\ref{solution}). Then, for all $t\geq0$,
\[
\mu^{(\lambda)}_{\infty,t}=\exp(-\lambda^{\alpha}t)\mu^{(\lambda
)}_{\infty}=m(t)\mu^{(\lambda)}_{\infty}.
\]
Reciprocally,\vspace*{1pt} if $(\mu_t, t \geq0)$ is a quasi-stationary solution of
the fragmentation equation, then there exists a $\lambda>0$ such that
$(\mu_t, t \geq0)=(\mu^{(\lambda)}_{\infty,t}, t \geq0)$.
\end{theorem}

\subsubsection*{Organization of the paper}
In Section \ref{background}, we begin with some background on subordinators and
regular variation. Section \ref{SectionSub} is the core of this paper:
our main results on large time behavior of self-similar Markov
processes conditioned on nonextinction are stated and proved there.
Together with Theorem \ref{theu}, these results imply Theorems \ref
{theq}, \ref{notcompact} and \ref{stat}, as well as Propositions \ref
{mproperties} and \ref{powernc}. Section \ref{SectionFS} is devoted
to the proof of Proposition \ref{fasterslower}.
Some properties of the limit measure $\mu_{\infty}$ are given in
Section \ref{PropertiesZ} and used to prove Proposition \ref{mass1}.
Finally, some specific examples are discussed in Section \ref
{Examples} and the proof of Theorem \ref{theu} is given in the \hyperref[App]{Appendix}.

\section{Background on subordinators and regular variation}
\label{background}

\subsection{Subordinators}
\label{sub}
A \textit{subordinator} is a nondecreasing L\'evy process, that is, a
nondecreasing c\`adl\`ag process with stationary and independent
increments. We recall here the main properties we need in this paper
and refer to Chapter 3 of \cite{BertoinLevy} for a more complete
introduction to the subject.

The distribution of a subordinator $(\xi_t, t \geq0)$ starting from
$\xi_0=0$ is characterized by its so-called \textit{Laplace exponent}
$\phi\dvtx[0,\infty[\ \rightarrow[0,\infty[$ via the identity
\[
\mathbb E [\exp(-\lambda\xi_t) ]=\exp(-t\phi(\lambda
))\qquad \forall\lambda,t \geq0.
\]
According to the L\'evy--Khintchine formula \cite{BertoinLevy},
Theorem 1, Chapter 1, there exists a real number $d\geq0$ and a measure
$\Pi$ on $]0,\infty[$, $\int_0^{\infty} (1 \wedge x) \Pi( d
x)<\infty$ such that
\[
\phi(\lambda)=d\lambda+ \int_0^{\infty}\bigl(1-\exp(-\lambda x)\bigr)\Pi
( d x)\qquad
\forall\lambda\geq0.
\]
The measure $\Pi$ governs the jumps of the subordinator: the jumps
process of $\xi$ is a Poisson point process with intensity $\Pi$.

We will need the strong Markov property of subordinators (\cite
{BertoinLevy}, Proposition~6, Chapter 1): given a subordinator $\xi$
and a stopping time $T$ with respect to the filtration $(\mathcaligr F_t,
t \geq0)$ generated by $\xi$, then, conditionally on $\{ T<\infty\}
$, the process $(\xi_{t+T} -\xi_T, t \geq0)$ is independent of
$\mathcaligr F_T$ and is distributed as $\xi$. Finally, we recall that
the semigroup of a subordinator possesses the Feller property
(\cite{BertoinLevy}, Proposition~5, Chapter 1).

\textit{Hereafter, all subordinators considered in this paper start
from $0$ and have drift $d=0$}. Their distribution is therefore
completely determined by their L\'evy measure $\Pi$. Note that when
$\Pi$ is related to a measure $B$ on $]0,1[$ via the formula (\ref
{PiB}), the above expression for $\phi$ coincides with that given by
equation (\ref{defphi}), that is,
\[
\phi(\lambda)=\int_0^{\infty}\bigl(1-\exp(-\lambda x)\bigr)\Pi( d
x)=\int_0^1(1-x^{\lambda})x B( dx)\qquad \forall\lambda\geq0.
\]
%
\subsection{Regular variation}
\label{secRV}
A function $f \dvtx]0,\infty[\  \rightarrow\,]0,\infty[$ is said to vary
regularly at $\infty$ (resp., 0) with index $\gamma\in\mathbb R$ if,
for all $a>0$,
\[
\frac{f(ax)}{f(x)} \rightarrow a^{\gamma} \qquad\mbox{ as }x \rightarrow
\infty\mbox{ (resp.,  0)}.
\]
We refer to \cite{BGT} for background on this topic. In particular, we
have already implicitly used the fact that the inverse, when it exists,
of a function regularly varying at $\infty$ with index $\gamma>0$ is
also regularly varying at $\infty$, with index $1/\gamma$ (see
Section~1.5.7 of \cite{BGT}).

Note that when the L\'evy measure $\Pi$ is related to the
fragmentation measure $B$ by the formula (\ref{PiB}), our main
assumption (\ref{hypothese}) reads
``$u \in\,]0,\infty[ \ \rightarrow\int_{u}^{\infty} \Pi( dx) $
varies regularly at $0$ with
index $-\beta$.'' It is classical that this is equivalent to the fact that
\[
\mbox{the function $\phi$ varies regularly at $\infty$ with index
$\beta$.}
\]
This can be easily proven using the Karamata Abelian--Tauberian theorems
(see, in particular, Chapters 1.6 and 1.7 of \cite{BGT}). We will
often use this form of the assumption (\ref{hypothese}).

To prove Theorem \ref{ThCVSub} below, which will then imply Theorems
\ref{theq} and \ref{notcompact}(ii), we will need the following
technical lemma, which is taken from Chow and Cuzick~\cite{ChowCuzick}.
\begin{lemma}[(Chow and Cuzick \cite{ChowCuzick}, Lemma 3)]
\label{CC}  Let $f$ be
regularly varying
at infinity with index $\gamma>0$ and suppose that for all
$\varepsilon>0$,
there exists some $x(\varepsilon)$ such that
%
%
\setcounter{equation}{10}
\renewcommand{\theequation}{\arabic{equation}}
\begin{equation}
\label{normalized}
\lambda^{\gamma-\varepsilon} \leq\frac{f(\lambda x)}{f(x)} \leq
\lambda^{\gamma+\varepsilon}\qquad \forall\lambda\geq1, \forall x \geq
x(\varepsilon).
\end{equation}
Then, for all $\theta>-1$,
\[
e^{f(t)} \biggl(\frac{f(t)}{t} \biggr)^{\theta+ 1} \int_t^{\infty}
(x-t )^{\theta}e^{-f(x)}\, d x
\mathop{\rightarrow}_{t \rightarrow\infty}
\gamma^{-1-\theta}\Gamma(1+\theta).
\]
\end{lemma}

We point out that Chow and Cuzick state their result for all regularly
varying functions with a positive index, but that their proof strongly
relies on the key point~(\ref{normalized}), which is not true for any
regularly varying function (counterexamples can easily be constructed).
However, the functions we are interested in, that is, $-\ln(m)$, and
to which we will apply this result, will, in general, satisfy (\ref
{normalized}). In particular, see Lemma \ref{lemmanormalized} below.

\section{Asymptotic behavior of self-similar Markov processes}
\label{SectionSub}

Given the construction (\ref{solution}) via subordinators of solutions
of the fragmentation equation, the issue of characterizing the large
time asymptotics of these solutions is equivalent to characterizing
large time behavior of distributions of time-changed subordinators.

So, let $\xi$ be a subordinator started from $0$ with L\'evy
measure $\Pi$ and no drift. We denote by $\phi$ its Laplace exponent.
Now, consider $\alpha<0$ and let $X(0)$ be a strictly positive random
variable, independent of $\xi$. Our goal is to specify the
asymptotic behavior as $t \rightarrow\infty$ of the distributions of
the random variables
%
%
\setcounter{equation}{11}
\renewcommand{\theequation}{\arabic{equation}}
\begin{equation}
\label{defX}
X(t):=X(0)\exp\bigl(-\xi_{\rho(X(0)^{\alpha}t )} \bigr),
\end{equation}
conditional on $\{X(t)>0\}$, where $\rho$ is given by
\[
\rho(t)=\inf\biggl\{ u \geq0 \dvtx\int_0^u \exp(\alpha\xi_r)
\, d r >t \biggr\}.
\]
Following Lamperti \cite{Lamperticaract}, the process $X$ belongs to
the so-called family of \textit{self-similar Markov processes}. This
means that it is strongly Markovian and that for all $x>0$, if $\mathbb
P_x$ denotes the distribution of $X$ started from $x$, then, for all $a>0$,
\[
\mbox{ the distribution of } \bigl(aX(a^{\alpha}t), t \geq0
\bigr) \mbox{ under } \mathbb P_x \mbox{ is } \mathbb P_{ax}.
\]
Moreover, $X$ reaches $0$ a.s. and it does so continuously. Conversely,
Lamperti \cite{Lamperticaract} also shows that any nonincreasing
c\`adl\`ag self-similar Markov processes on $[0,\infty[$ that reaches 0
continuously in finite time a.s. can be constructed in this way via a
time-changed subordinator.

Note that the moment at which $X$ reaches $0$ is $X(0)^{|\alpha|}I$,
where $I$ is the \textit{exponential functional} defined by
%
%
\begin{equation}
\label{defI}
I:=\int_0^{\infty} \exp(\alpha\xi_r) \, d r,
\end{equation}
which is clearly a.s. finite. The distribution of the random variable
$I$ was first studied in detail in
\cite{CarmonaPetitYor}. In particular, it is well known that for all
integers $n \geq1$,
\[
\mathbb E[I^n]=\frac{n !}{\phi(|\alpha|)\phi(2|\alpha|)\cdots\phi(n
|\alpha|)},
\]
and that the distribution of $I$ is characterized by these moments
(\cite{CarmonaPetitYor}, Proposition~3.3). It will also be essential
for us (see \cite{BYFacExp}, Propositions 1 and 2) that there exists a unique
probability measure $\mu_R$ on $]0,\infty[$ whose entire positive
moments are
given by
%
%
\begin{equation}
\label{defmuR}
\int_0^{\infty} x^n \mu_R ( d
x)=\phi(|\alpha|)\phi(2|\alpha|)\cdots\phi(n|\alpha|),\qquad n \geq1,
\end{equation}
and that, moreover, if $R$ denotes a random variable with distribution
$\mu_R$
independent of $I$, then
%
%
\begin{equation}
\label{face}
RI \stackrel{\dd}{=}\mathbf e (1),
\end{equation}
where $\mathbf e(1)$ has an exponential
distribution with parameter $1$.

We now have the material necessary to state the main result of this section.
To be consistent with the notation used for the fragmentation equation,
we denote by $x\mu_0( d x)$, $x>0$, the distribution of $X(0)$.
Also, we recall the definition of the function $\varphi$ as the
inverse, well defined in a neighborhood of $\infty$, of $t \rightarrow
t/\phi(t)$.
\begin{theorem}
\label{ThCVSub} Suppose that $\int_{u}^{\infty} \Pi( dx) $
varies regularly at $0$ with index $-\beta$, $\beta\in[0,1[$, and
$\int^{\infty} x \Pi( d x)<\infty$.
\begin{longlist}
\item If the support of $\mu_0$ is bounded with a supremum
equal to $1$, then, for all
bounded continuous functions $f\dvtx]0,\infty[ \ \rightarrow\mathbb R$,
\[
\mathbb
E \biggl[f \biggl( \biggl(\frac{\varphi(|\alpha|t)}{|\alpha|t}
\biggr)^{1/|\alpha|}X(t) \biggr)
\Big| X(t)>0 \biggr] \mathop{\rightarrow}_{t \rightarrow\infty}
\mathbb E \bigl[
f\bigl(R^{1/|\alpha|}\bigr) \bigr],
\]
where $R$ is the random variable with distribution $\mu_R$ defined by
(\ref{defmuR}).
\item
If $\mu_0$ possesses a density $u_0$ in a neighborhood
of $\infty$ such that
\[
\ln( u_0(x) ) \mathop{\sim}_{\infty}-Cx^{\gamma}
\]
for some $\gamma>0$, then the function $t \in\,]0,\infty[\  \rightarrow
\mathbb P(X(t)>0)$ is continuously differentiable. If, moreover, the
derivative of $t \rightarrow\ln(\mathbb P(X(t)>0))$ is regularly
varying at $\infty$---which is true when, for example, this derivative
is monotone near $\infty$---then, for all
bounded continuous functions $f\dvtx]0,\infty[\  \rightarrow\mathbb R$, as
$t \rightarrow\infty$,
\[
\mathbb
E \biggl[f \biggl( \biggl(\frac{h(t)}{C_{\alpha,\beta,\gamma
,C}t} \biggr)^{1/|\alpha|}X(t) \biggr) \Big| X(t)>0
\biggr] \mathop{\rightarrow}_{t \rightarrow\infty} \mathbb E \bigl[
f\bigl(R^{1/|\alpha|}\bigr) \bigr],
\]
where the function $h$ is the inverse, defined in the neighborhood of
$\infty$, of $t \rightarrow t^{1+|\alpha|/\gamma}/\phi(t)$ and
$C_{\alpha,\beta,\gamma,C}$ is the constant defined in the statement
of Theorem~\ref{notcompact}.
\end{longlist}
\end{theorem}

We will see in the proof of this result that the function $t
\rightarrow\varphi(|\alpha|t)/|\alpha|t$ in assertion (i) can be
replaced by any asymptotically equivalent function and likewise for $h$
in the second assertion.

Now, let $B$ be the fragmentation measure related to $\Pi$ by (\ref
{PiB}). If $(\mu_t , t \geq0)$ refers to the solution of the $(\alpha
,B)$-fragmentation equation constructed from $\xi$ by the formula
(\ref{solution}), then we have
\[
m(t)=\int_0^{\infty} x \mu_t( dx)=\mathbb P \bigl(
X(t)>0 \bigr)
\]
and the distribution of $X(t)$ conditional on $X(t)>0$ is $x\mu
_t( dx)/m(t)$. The above theorem then leads directly to the
statements of Theorems \ref{theq} and \ref{notcompact}(ii)
[note that $\int^{\infty}x\Pi( dx)<\infty$ is equivalent to
$\int_0 |{\ln(x)}|x B( dx)<\infty$]. The limit distribution
$x\mu_{\infty}( dx)$ mentioned in these theorems is therefore
the distribution of $R^{1/|\alpha|}$. The large time behavior of
$m(t)=\mathbb P(X(t)>0)$ is studied in Section \ref{totalmass} below,
whereas Theorem \ref{ThCVSub} is established in Section \ref{ProofCVSub}.

We finish with the following result on the \textit{quasi-stationary
distributions} of $X$, which will be proven in Section \ref{QS} and
which, in terms of the fragmentation equation, will lead to Theorem
\ref{stat}. We recall that the quasi-stationary distributions of $X$
are the distributions
$\varsigma$ on $]0,\infty[$ such that
\begin{eqnarray}
X(0) \stackrel{\dd}{\sim} \varsigma\quad\Rightarrow\quad\mathbb
E[f(X(t)) | X(t)>0]=\mathbb
E[f(X(0))]\nonumber
\end{eqnarray}
for all $t \geq0$ and all test functions $f$ defined on
$]0,\infty[$.

\begin{theorem}
\label{ThQS}
Let $\mu_R^{(\lambda)}$ denote the law of $\lambda
R^{1/|\alpha|}$, $\lambda>0$. Then, a probability measure $\varsigma
$ on $]0,\infty[$ is a quasi-stationary
distribution of $X$ if and only if $\varsigma=\mu_R^{(\lambda)}$ for
some $\lambda>0$. Moreover, if $X(0) \stackrel{\dd}{\sim} \mu
_R^{(\lambda)}$, then
\[
\mathbb P\bigl(X(t)>0\bigr)=\exp(-\lambda^{\alpha} t)\qquad \forall t \geq0.
\]
\end{theorem}

We point out that this theorem does not lead directly to the reciprocal
assertion of Theorem \ref{stat}. However, easy manipulations of the
fragmentation equation will lead to it; see Section \ref{QS} for details.

\subsection{Total mass behavior}
\label{totalmass}
This section is devoted to the description of the behavior of the total mass
\[
m(t)=\int_0^{\infty} x \mu_t( dx)=\mathbb P \bigl(
X(t)>0 \bigr)=\mathbb P \bigl(I >X(0)^{\alpha}t \bigr).
\]
The notation is that introduced above in the introduction of Section
\ref{SectionSub}. We start with the following result, which holds for
all fragmentation equations with parameters $\alpha<0, B$ and all
initial measures $\mu_0$ such that $\int_0^{\infty}x \mu_0( dx)=1$.
\begin{proposition}
\label{massgeneral}
The total mass $m$ is strictly positive and strictly decreasing on
$[0,\infty[$. Moreover, $m(t) \rightarrow0$ as $t \rightarrow\infty$.
\end{proposition}
\begin{pf}
Since
\[
m(t)=\int_0^{\infty} \mathbb P(I>x^{\alpha} t) x\mu_0( dx),
\]
it is sufficient to show that the function $t \in[0,\infty[\  \rightarrow\mathbb P(I>t)$ is strictly positive, strictly decreasing
and converges to $0$ as $t \rightarrow\infty$. This last point is
obvious since $I<\infty$ a.s. Next, suppose that $\mathbb P(I \leq
t)=1$ for some $t>0$. This would imply that for all $n \geq1$,
\[
\frac{n !}{\phi(|\alpha|)\phi(2|\alpha|)\cdots\phi(n |\alpha
|)}=\mathbb E[I^n] \leq t^n.
\]
However, we saw in the \hyperref[sec1]{Introduction} that $x/\phi(x) \rightarrow\infty
$ as $x \rightarrow\infty$. In particular, $2t \leq n/\phi(n |\alpha
|)$ for large enough $n$, say $n> n_0$. Hence, we would have
\[
\frac{n_0 !}{\phi(|\alpha|)\phi(2|\alpha|)\cdots\phi(n_0 |\alpha
|)}(2t)^{n-n0} \leq t^n
\]
for all $n>n_0$, which is impossible. Therefore,
$\mathbb P(I>t)>0$ for all $t>0$.

Finally, for all $t>0$, using the Markov property of subordinators, we get
\begin{eqnarray*}
I &= &\int_0^t \exp(\alpha\xi_r) \, dr+\exp(\alpha\xi
_t)\int_0^{\infty} \exp\bigl(\alpha(\xi_{r+t}-\xi_t)\bigr) \, dr \\
& \leq& t+ \exp(\alpha\xi_t) \tilde I,
\end{eqnarray*}
where $\tilde I$ is distributed as $I$ and is independent of $\xi
_{t}$. Consider $a$ such that $\mathbb P(I \leq a)>0$ and note, using
the Poisson point process construction of the subordinator, that
$\mathbb P(\exp(\alpha\xi_t) \leq t/a)>0$ for all $t>0$. Then,
\[
0<\mathbb P \bigl( \exp(\alpha\xi_t) \leq t/a, \tilde I \leq a
\bigr) \leq\mathbb P(I \leq2t)\qquad \forall t >0.
\]
This leads to the fact that $\mathbb P(t \geq I>s)>0$ for all $0 \leq
s<t$. Indeed, the event $\{I>s \}$ coincides with $\{\rho(s)<\infty\}
$ and when $I>s$,
\[
I=s+\exp\bigl(\alpha\xi_{\rho(s)}\bigr)\int_0^{\infty}\exp\bigl(\alpha
\bigl(\xi_{r+\rho(s)}-\xi_{\rho(s)} \bigr)\bigr) \, dr.
\]
Using the strong Markov property of the subordinator at the stopping
time $\rho(s)$, we get, with probability 1,
\[
(I-s)^+=\exp\bigl(\alpha\xi_{\rho(s)}\bigr) \tilde I
\]
with $\tilde I$ independent of $\xi_{\rho(s)}$ and distributed as
$I$. Hence, for all $0 \leq s <t$,
\begin{eqnarray*}
\mathbb P(I>s) -\mathbb P(I>t)&=&\mathbb P(s<I \leq t)\\
&=& \mathbb P
\bigl(\exp\bigl(\alpha\xi_{\rho(s)}\bigr)>0, \tilde I \leq(t-s)\exp\bigl(|\alpha|
\xi_{\rho(s)}\bigr) \bigr)
\end{eqnarray*}
and this last probability is strictly positive since $\mathbb P(\exp
(\alpha\xi_{\rho(s)})>0)=\mathbb P(I>s)>0$ and $\mathbb P(\tilde I
\leq a)>0$ for all $a>0$.
\end{pf}

We now turn to the proofs of the more precise descriptions of the
behavior of $m$ stated in Proposition \ref{mproperties}, Theorem \ref
{notcompact}(i) and Proposition \ref{powernc}.
The crucial point is the following lemma, which is basically a
consequence of Rivero \cite{VictorLog}, Proposition~2, and K\"onig and
M\"orters \cite{KM}, Lemma 2.3.
\begin{lemma}
\label{Victor}
Assume that (\ref{hypothese}) holds
or, equivalently, that
$\phi$ varies regularly at $\infty$ with index $\beta\in[0,1[$.
Then,
\[
-\ln\bigl(\mathbb P(I>t) \bigr) \mathop{\sim}_{\infty} \frac
{(1-\beta)}{|\alpha|}
\varphi(|\alpha|t) \mathop{\sim}_{\infty} (1-\beta)|\alpha
|^{{\beta}/({1-\beta})} \varphi(t),
\]
where $\varphi$ is the inverse of $t \rightarrow
t/\phi(t)$, which is well defined in the neighborhood of~$\infty$.
In particular, $-\ln(\mathbb P(I>t) )$ is regularly
varying at $\infty$ with index $1/(1-\beta)$.
\end{lemma}
\begin{pf}
Note that the Laplace exponent of the
subordinator $|\alpha|\xi$ is $\phi(|\alpha| \cdot)$ and that the
inverse of $t \rightarrow t/\phi(|\alpha| t)$ is $\varphi(|\alpha|
\cdot)/|\alpha|$. Using these facts, we can restrict our proof to the
case $|\alpha|=1$, which is supposed in the following.

When $\beta\in\,]0,1[$, the statement of the lemma is exactly
Proposition 2 of Rivero~\cite{VictorLog}. When $\beta=0$ and $\phi
(\infty)<\infty$,
\[
\frac{1}{n}\ln\biggl(\frac{ \mathbb E(I^n)}{n!} \biggr)=-\frac
{1}{n}\sum^n_{i=1}\ln(\phi(i))\xrightarrow{n\to\infty} -\ln\phi
(\infty).
\]
Then, by Lemma 2.3. of K\"onig and M\"orters \cite{KM},
\[
\lim_{t\to\infty}\frac{1}{t}\ln\bigl( \mathbb P(I>t)
\bigr)=-\phi(\infty).
\]
Finally, when $\beta=0$ and $\phi(\infty)=\infty$, we can adapt K\"
onig and M\"orters' proof of~\cite{KM}, Lemma 2.3, to obtain the
expected result. Indeed, first note that
%
%
\begin{equation}\qquad
\label{lim1}
\frac{1}{n}\ln\biggl( \mathbb E \biggl[ \frac{I^n \phi(n)
^n}{n^n} \biggr] \biggr)=\frac{1}{n}\ln\biggl(\frac{n!}{n^n}
\biggr) + \ln(\phi(n))-\frac{1}{n} \sum_{i=1}^{n} \ln(\phi(i)) \mathop
{\rightarrow}_{n \rightarrow\infty} -1
\end{equation}
as a consequence of Stirling's formula and of the fact that
\[
\ln(\phi(n))-\frac{1}{n} \sum_{i=1}^{n} \ln(\phi(i))\mathop
{\rightarrow}_{n \rightarrow\infty} 0
\]
since $\phi$ is a slowly varying function (see Section 3.2 of Rivero
\cite{VictorLog} for a proof of this last point). It is then easy,
using Markov's inequality, to show that
\[
\limsup_{n \rightarrow\infty} \frac{1}{n}\ln\bigl( \mathbb P
\bigl( I >n/\phi(n) \bigr) \bigr) \leq-1.
\]
To get a lower bound for the limit inferior, set $Y_n:=\ln(I\phi
(n)/n)$. For every $\varepsilon>0$ and every integer $m$, we have that
\[
\frac{1}{\mathbb E[I^n]} \mathbb E \bigl[ I^n \mathbf1_{\{Y_n\geq
\varepsilon\}} \bigr] \leq\exp(-\varepsilon m)
\frac{\mathbb E[I^{n+m}] \phi(n)^m}{\mathbb E[I^n] n^m}
\mathop{\rightarrow}_{n \rightarrow\infty}\exp(-\varepsilon m) .
\]
Letting $m \rightarrow\infty$, this gives
%
%
\begin{equation}
\label{lim2}
\frac{1}{\mathbb E[I^n]} \mathbb E \bigl[ I^n \mathbf1_{\{Y_n\geq
\varepsilon\}} \bigr] \mathop{\rightarrow}_{n \rightarrow\infty} 0.
\end{equation}
Besides, for all $\varepsilon>0$ and all $n \geq1$,
\[
\frac{1}{n} \ln\bigl( \mathbb P\bigl(I>n\exp(-\varepsilon)/ \phi(n)\bigr)\bigr)
\geq\frac{1}{n} \ln\bigl( \mathbb P(|Y_n| < \varepsilon)\bigr).
\]
However, $I^{-n} > \exp(-n\varepsilon) n^{-n}\phi(n)^n$ on $\{ |Y_n|
< \varepsilon\}$, which gives
\begin{eqnarray*}
\frac{1}{n} \ln\bigl( \mathbb P(|Y_n| < \varepsilon) \bigr)&=&
\frac{1}{n} \ln\biggl(\frac{\mathbb E[I^{-n}I^n \mathbf1_{\{ |Y_n|
< \varepsilon\}} ]}{\mathbb E[I^n]}\mathbb E[I^n] \biggr)\\
& \geq& \frac{1}{n} \ln\biggl( \exp(-n\varepsilon) \frac{\mathbb
E[I^n \mathbf1_{\{ |Y_n| < \varepsilon\}} ]}{\mathbb E[I^n]}
n^{-n}\phi(n)^n\mathbb E[I^n] \biggr).
\end{eqnarray*}
By (\ref{lim1}) and (\ref{lim2}), the last line of this inequality
converges to $-\varepsilon-1$ as $n \rightarrow\infty$.
Thus, since the function $t \rightarrow t/\phi(t)$ is increasing and
$\varphi(t)\rightarrow\infty$ as $t \rightarrow\infty$, we have
proven that
\begin{eqnarray*}
\limsup_{t \rightarrow\infty} \frac{1}{\varphi(t)}\ln\bigl(
\mathbb P ( I >t ) \bigr)&\leq&-1,
\\
\liminf_{t \rightarrow\infty} \frac{1}{\varphi(t\exp(\varepsilon
))}\ln\bigl( \mathbb P ( I >t ) \bigr) &\geq&
-\varepsilon-1.
\end{eqnarray*}
Using the regular variation of $\varphi$, we get the expected
\[
\lim_{t\to\infty}\frac{1}{\varphi(t)}\ln\bigl( \mathbb
P(I>t) \bigr)=-1.
\]
\upqed\end{pf}

\subsubsection{$\mu_0$ with bounded support: Proof of Proposition
\protect\ref{mproperties}}
\label{subsectionproofcomp}

We recall that, with no loss of generality, the supremum of the support
of $\mu_0$ is supposed to be equal to $1$.
Thus,
\[
m(t)=\int_0^{1} \mathbb P ( I >t x^{\alpha} ) x \mu_0
( d x) \leq\mathbb P ( I>t ) \int_0^{1} x \mu_0
( d x) = \mathbb P ( I>t ).
\]
According to Proposition 3.3 in \cite{CarmonaPetitYor}, $C_{\lambda
}:=\mathbb E[\exp(\lambda I)]<\infty$ provided that $\lambda<\phi
(\infty)$. Hence, for such $\lambda$'s,
\[
m(t) \leq\mathbb P ( I>t ) \leq C_{\lambda}\exp
(-\lambda t)\qquad \forall t \geq0,
\]
which gives the first part of the statement.

Now, assume that (\ref{hypothese}) holds. Then, on the one hand, since
$m(t) \leq\mathbb P(I>t)$, we get,
by Lemma \ref{Victor},
\[
\liminf_{t \rightarrow\infty} \frac{-\ln( m(t)
)}{\varphi(|\alpha|t)} \geq\frac{1-\beta}{|\alpha|}.
\]
On the other hand, for all $0<\varepsilon<1$,
\[
m(t) \geq\mathbb P \bigl( I >t (1-\varepsilon)^{\alpha} \bigr)
\int_{1-\varepsilon}^{1} x \mu_0 ( d x).
\]
By assumption, $\int_{1-\varepsilon}^{1} x \mu_0 ( d x)>0$, hence
\[
\limsup_{t \rightarrow\infty} \frac{-\ln( m(t) )
}{\varphi(|\alpha|t)} \leq\limsup_{t \rightarrow\infty} \frac
{-\ln( \mathbb P(I>t(1-\varepsilon)^{\alpha})
)}{\varphi(|\alpha| t)} =\frac{1-\beta}{|\alpha|}
(1-\varepsilon)^{{\alpha}/({1-\beta})}.
\]
Then, let $\varepsilon\downarrow0$ to get the expected result.

\subsubsection{$\mu_0$ with unbounded support: Proofs of Theorem
\protect\ref{notcompact}\textup{(i)} and Proposition~\protect\ref{powernc}}

\mbox{}

\begin{pf*}{Proof of Theorem \protect\ref{notcompact}\textup{(i)}}
First, suppose that $\mu_0( d x)= \exp(-C x^{\gamma}) \, d x$,
$\gamma>0$. We have
%
%
\begin{eqnarray}
\label{regdev}
m(t)&=&\int_0^{\infty} \mathbb P ( I >t x^{\alpha} ) x
\exp(-Cx^{\gamma}) \, d x\nonumber\\[-8pt]\\[-8pt]
&=& \frac{t^{-2/\alpha}}{\gamma}\int
_0^{\infty} \mathbb P ( I >u^{\alpha/\gamma} )
u^{2/\gamma-1} \exp( -Cut^{-\gamma/\alpha} ) \, d u,\nonumber
\end{eqnarray}
using the change of variable $u=(xt^{1/\alpha})^{\gamma}$. Now, use
Lemma \ref{Victor} and Theorem~4.12.10(iii) of \cite{BGT} to get
\begin{eqnarray*}
-\ln\biggl(\int_0^x \mathbb P ( I >u^{\alpha/\gamma} )
u^{2/\gamma-1} \, d u \biggr) &\mathop{\sim}\limits_{x\rightarrow0}&
-\ln\bigl( \mathbb P ( I>x^{\alpha/\gamma} ) \bigr)\\
&\mathop{\sim}\limits_{x\rightarrow0}& (1-\beta)|\alpha|^{\beta/{(1-\beta)}}
\varphi( x^{\alpha/\gamma}),
\end{eqnarray*}
which varies regularly at $0$ with index $\alpha/{(\gamma(1-\beta
))}$. Note that in a neighborhood of $0$, $x \rightarrow1/\varphi(
x^{\alpha/\gamma})$ is the inverse of
\[
x \rightarrow\biggl(x \phi\biggl(\frac{1}{x} \biggr)
\biggr)^{-\gamma/\alpha}.
\]
Hence, by de Bruijn's Tauberian theorem (\cite{BGT}, Theorem 4.12.9) we have
\[
-\ln\biggl( \int_0^{\infty} \mathbb P ( I >u^{\alpha/\gamma}
) u^{2/\gamma-1} \exp( -ut ) \, d u \biggr)
\mathop{\sim}_{t \rightarrow\infty}C_{\alpha,\beta,\gamma} /h_0 (t),
\]
where $h_0$ is the inverse, well defined in the neighborhood of $\infty
$, of $x \rightarrow x^{-1} (x \phi(1/x) )^{\gamma/\alpha
}$ and
$
C_{\alpha,\beta,\gamma}$ is the constant defined in the statement of
Theorem~\ref{notcompact}(i). Together with (\ref{regdev}), this leads to
\[
-\ln( m(t) ) \mathop{\sim}_{\infty} C_{\alpha,\beta
,\gamma} C^{(1+(1-\beta)\gamma/|\alpha|)^{-1}} / h_0\bigl(t^{\gamma
/|\alpha|}\bigr).
\]
In other words,
\[
-\ln( m(t) ) \mathop{\sim}_{\infty} C_{\alpha,\beta
,\gamma} C^{(1+(1-\beta)\gamma/|\alpha|)^{-1}} h (t),
\]
where $h$ is the inverse of $t^{1+|\alpha|/\gamma}/\phi(t)$.

Now, suppose that $\mu_0$ possesses a density $u_0$ in a
neighborhood of $\infty$ such that $\ln( u_0(x) )
\mathop{\sim}_{\infty}-Cx^{\gamma}$, $\gamma>0$. Fix $\varepsilon
>0$ and let $C_{\varepsilon}$ be such that $u_0(x)$ exists for $x \geq
C_{\varepsilon}$ and
%
%
\begin{equation}
\label{encadr}
\exp\bigl(-(1+\varepsilon)Cx^{\gamma}\bigr)
\leq u_0(x) \leq\exp\bigl(-(1-\varepsilon)Cx^{\gamma}\bigr) \qquad \forall x \geq
C_\varepsilon.
\end{equation}
Then, write
\[
m(t)=\int_0^{C_{\varepsilon}} \mathbb P ( I >t x^{\alpha}
) x \mu_0( dx) +\int_{C_{\varepsilon}}^{\infty}
\mathbb P ( I >t x^{\alpha} ) x u_0(x) \, dx.
\]
On the one hand, following the argument developed in Section \ref
{subsectionproofcomp}, we get
\[
\limsup_{t \rightarrow\infty} \frac{\ln( \int
_0^{C_{\varepsilon}} \mathbb P ( I >t x^{\alpha} ) x \mu
_0( dx) )}{\varphi(t)} \leq-(1-\beta)|\alpha|^{\beta
/(1-\beta)}C_{\varepsilon}^{\alpha/(1-\beta)},
\]
which actually holds for any initial measure $\mu_0$. Note that
$\varphi(t)/h(t) \rightarrow\infty$ as $t \rightarrow\infty$,
where $h$ is the function defined above in the first part of this proof.

On the other hand, inequalities (\ref{encadr}) and the results of
the first part of this proof imply that
\[
\limsup_{t \rightarrow\infty} \frac{-\ln( \int
_{C_{\varepsilon}}^{\infty} \mathbb P ( I >t x^{\alpha}
) x u_0(x) \, dx ) }{h(t)} \leq C_{\alpha,\beta,\gamma}
\bigl((1+\varepsilon)C\bigr)^{(1+(1-\beta)\gamma/|\alpha|)^{-1}}
\]
and
\[
\liminf_{t \rightarrow\infty} \frac{-\ln( \int
_{C_{\varepsilon}}^{\infty} \mathbb P ( I >t x^{\alpha}
) x u_0(x) \, dx ) }{h(t)} \geq C_{\alpha,\beta,\gamma}
\bigl((1-\varepsilon)C\bigr)^{(1+(1-\beta)\gamma/|\alpha|)^{-1}}.
\]
We have therefore proven that
\begin{eqnarray*}
C_{\alpha,\beta,\gamma} \bigl((1-\varepsilon)C\bigr)^{(1+(1-\beta)\gamma
/|\alpha|)^{-1}} &\leq&\liminf_{t \rightarrow\infty} \frac{-\ln
(m(t)) }{h(t)} \leq\limsup_{t \rightarrow\infty} \frac{-\ln
(m(t))} {h(t)} \\
&\leq& C_{\alpha,\beta,\gamma} \bigl((1+\varepsilon
)C\bigr)^{(1+(1-\beta)\gamma/|\alpha|)^{-1}}
\end{eqnarray*}
for all $\varepsilon>0$. The result follows by letting $\varepsilon
\downarrow0$.
\end{pf*}
\begin{pf*}{Proof of Proposition \protect\ref{powernc}}
Suppose that $u_0(x)=Cx^{-\gamma}$ on
$[a,\infty[$ for some $a>0$ and $\gamma>2$. Then,
\[
m(t)= C\int_a^{\infty} \mathbb P(I>x^{\alpha}t) x^{1-\gamma}
\, d x + \int_0^{a} \mathbb P(I>x^{\alpha}t) x\mu_0( dx).
\]
With the change of variables $u=x^{\alpha}t$,
\[
\int_a^{\infty} \mathbb P(I>x^{\alpha}t) x^{1-\gamma} \, d
x=\frac{t^{({\gamma-2})/{\alpha}} }{|\alpha|}\int_0^{a^{\alpha
}t} \mathbb P(I>u) u^{({2-\gamma})/{\alpha}-1} \, du
\]
and this last integral converges to a finite limit as $t \rightarrow
\infty$ since $\mathbb P(I>u) \leq C_{\lambda} \exp(-\lambda u)$ for
all $u \geq0$ and some $\lambda>0$ sufficiently small (see the proof
of Proposition \ref{mproperties} for this last point). Using the same
upper bound for $\mathbb P(I >x^{\alpha} t)$, we get that
\[
\int_0^{a} \mathbb P(I>x^{\alpha}t) x\mu_0( dx) \leq
C_{\lambda} \exp(-\lambda a^{\alpha}t) \int_0^a x \mu_0( dx).
\]
Thus,
\[
m(t) \mathop{\sim}_{t \rightarrow\infty} \frac{C}{|\alpha|}
t^{({\gamma-2})/{\alpha}} \int_0^{\infty} \mathbb P(I>u)
u^{({2-\gamma})/{\alpha}-1} \, du.
\]
It is not hard to extend this proof to the case where $u_0(x)\mathop
{\sim}_{\infty}Cx^{-\gamma} $, for some $\gamma>2$. This is left to
the reader.
\end{pf*}

\subsection{Proof of Theorem \protect\ref{ThCVSub}}
\label{ProofCVSub}

We start with the following lemma.
\begin{lemma}
\label{lemmaCVSub}
Suppose that $-\ln(m)$ varies regularly at $\infty$ with a positive
index $\gamma$ and satisfies (\ref{normalized}). Then, for any
function $g\dvtx[0,\infty[\ \rightarrow\,]0,\infty[$ such that $g(t) / (-\ln
(m(t))) \rightarrow1$ as $t \rightarrow\infty$, we have
\[
\mathbb
E \biggl[f \biggl( \biggl(\frac{\gamma g(t)}{t} \biggr)^{1/|\alpha
|}X(t) \biggr) \Big| X(t)>0
\biggr] \mathop{\rightarrow}_{t \rightarrow\infty} \mathbb E \bigl[
f\bigl(R^{1/|\alpha|}\bigr) \bigr]
\]
for all continuous bounded test functions $f$ on $]0,\infty[$.
\end{lemma}
\begin{pf}
First, note that when $X(0)^{|\alpha|} I>t$,
we have
\begin{eqnarray*}
X(0)^{|\alpha|} I
&=& X(0)^{|\alpha|} \int_0^{\rho(X(0)^{\alpha} t)}
\exp(\alpha\xi_r) \, d r\\
&&{} + X(0)^{|\alpha|} \exp\bigl(\alpha\xi
_{\rho(X(0)^{\alpha} t)}\bigr) \int
_0^{\infty}\exp\bigl(\alpha\bigl(\xi_{r+\rho(X(0)^{\alpha} t)}-\xi_{\rho
(X(0)^{\alpha} t)}\bigr)\bigr) \, d
r \\
&=& t+X(0)^{|\alpha|} \exp\bigl(\alpha\xi_{\rho(X(0)^{\alpha} t)}\bigr)
\int
_0^{\infty}\exp\bigl(\alpha\bigl(\xi_{r+\rho(X(0)^{\alpha} t)}-\xi_{\rho
(X(0)^{\alpha} t)}\bigr)\bigr) \, d
r.
\end{eqnarray*}
Now, use the strong Markov property of $\xi$ at the
(randomized) stopping time $\rho(X(0)^{\alpha} t)$ to get
%
%
\begin{equation}
\label{keyidentity}
\bigl(X(0)^{|\alpha|} I-t\bigr)^+=X(0)^{|\alpha|} \exp\bigl(\alpha\xi_{\rho
(X(0)^{\alpha} t)}\bigr) \tilde I=X(t)^{|\alpha|}\tilde I,
\end{equation}
where $\tilde I$ is distributed as $I$ and is
independent of $X(t)$.
This gives, for all $n
\in\mathbb N^{*}$,
\begin{eqnarray*}
&&
m(t)^{-1}\mathbb
E \biggl[ \biggl( \biggl(\frac{\gamma g(t)}{t} \biggr)^{1/|\alpha|}X(t)
\biggr)^{|\alpha| n} \biggr]\mathbb E[I^n]
\\
&&\qquad=
m(t)^{-1} \biggl(\frac{\gamma g(t)}{t} \biggr)^{n}\mathbb
E \bigl[X(t)^{|\alpha| n} \bigr]\mathbb E[I^n]\\
&&\qquad=m(t)^{-1} \biggl(\frac{\gamma g(t)}{t} \biggr)^{ n}\mathbb
E \bigl[ \bigl(\bigl(X(0)^{|\alpha|}I-t\bigr)^+ \bigr)^n \bigr].
\end{eqnarray*}
Then, recall that
\[
m(t)=\mathbb P\bigl(X(0)^{|\alpha|} I>t\bigr),\qquad t \geq0.
\]
Integrating by parts,
we have
\begin{eqnarray*}
&&
m(t)^{-1} \biggl(\frac{\gamma g(t)}{t} \biggr)^{ n}\mathbb
E \bigl[ \bigl(\bigl(X(0)^{|\alpha| }I-t\bigr)^+ \bigr)^n \bigr] \\
&&\qquad
= n m(t)^{-1} \biggl(
\frac{\gamma g(t)}{t} \biggr)^{n}\int_t^{\infty}(x-t)^{n-1} m(x)
\, d x,
\end{eqnarray*}
which, according to Lemma \ref{CC} and the assumptions we have made on
$-\ln(m)$ and $g$, converges as $t \rightarrow\infty$ to
$n!$.
Next, note that $\mathbb E[R^n]\mathbb E[I^n] = n !$, using
the factorization property (\ref{face}) of the exponential random
variable with
parameter $1$. Putting all of the pieces
together, we have proven that for all integers $n \geq1$,
\[
\mathbb
E \biggl[ \biggl( \biggl(\frac{\gamma g(t)}{t} \biggr)^{1/|\alpha|}X(t)
\biggr)^{|\alpha| n} \Big| X(t)>0 \biggr] \mathop{\rightarrow
}_{t \rightarrow\infty} \mathbb
E[R^n].
\]

\textit{Summary.} Let $\nu_t$ denote the distribution of $\gamma
t^{-1}g(t) X(t)^{|\alpha|}$ conditional on $X(t)>0$ ($\nu_t$ is a
probability measure on $]0,\infty[$). We have shown that for all $n
\geq1$,
\[
\int_{0}^{\infty}x^n \nu_t( dx) \rightarrow\int_{0}^{\infty
}x^n \mu_R( dx),
\]
where $\mu_R$ is the distribution of $R$. Of course, this still holds
for $n=0$, but the distribution of $R$ is characterized by its moments.
It is then well known (\cite{Feller}, Chapter VIII, page 269) that
this implies that $\nu_t$ converges in distribution to $\mu_R$.
\end{pf}

\subsubsection{Proof of Theorem \protect\ref{ThCVSub}\textup{(i)}}
\label{proofnormalized}
By Proposition \ref{mproperties}, under the hypothesis (\ref
{hypothese}), $-\ln(m)$ varies regularly at $\infty$ with index
$1/(1-\beta)$ and, more precisely,
\[
-\ln(m(t)) \mathop{\sim}_{t \rightarrow\infty} \frac{(1-\beta
)}{|\alpha|}
\varphi(|\alpha|t).
\]
Together with Lemma \ref{lemmaCVSub}, this implies the statement of
Theorem \ref{ThCVSub}, provided that $-\ln(m)$ satisfies (\ref
{normalized}). The goal of this section is to prove this last point
when $\mu_0$ has bounded support.
\begin{lemma}
\label{lemmanormalized}
Let
\[
f(x)=-\ln(m(x)),\qquad x\geq0,
\]
and assume that (\ref{hypothese}) holds, $\int^{\infty}x \Pi
( d x)<\infty$ and $\mu_0$ has bounded support. Then, for all
$\varepsilon>0$, there exists some $x(\varepsilon)$ such that
\[
\lambda^{{1}/({1-\beta})-\varepsilon}\leq\frac{f(\lambda
x)}{f(x)} \leq\lambda^{{1}/({1-\beta})+\varepsilon}\qquad
\forall\lambda\geq1 \mbox{ and } \forall x \geq x(\varepsilon).
\]
\end{lemma}

This lemma is a direct consequence of Lemmas \ref{tech1} and \ref
{tech4} below.
\begin{lemma}
\label{tech1} Let $g\dvtx]0,\infty[\  \rightarrow\,]0,\infty[$ be a
continuously differentiable function such that
\[
\frac{xg'(x)}{g(x)}\rightarrow c >0 \qquad\mbox{as } x \rightarrow\infty.
\]
Then, for all $\varepsilon>0$, there exists some $x(\varepsilon)$
such that
\[
\lambda^{c-\varepsilon}\leq\frac{g(\lambda x)}{g(x)} \leq\lambda
^{c+\varepsilon}\qquad
\forall\lambda\geq1 \mbox{ and } \forall x \geq x(\varepsilon).
\]
\end{lemma}
\begin{pf}
For $\varepsilon>0$, let $x(\varepsilon)$
be such that
\[
c-\varepsilon\leq\frac{xg'(x)}{g(x)} \leq c+\varepsilon\qquad\mbox{for
all } x \geq x(\varepsilon).
\]
For such $x$'s and all $\lambda\geq1$,
\begin{eqnarray*}
(c-\varepsilon) \ln(\lambda)
&=&(c-\varepsilon) \int_x^{\lambda x} y^{-1} \, d y \leq\int
_x^{\lambda x} \frac{g'(y)}{g(y)}\, d y\\
& \leq&(c+\varepsilon)
\int_x^{\lambda x} y^{-1} \, d y= (c+\varepsilon) \ln(\lambda).
\end{eqnarray*}
Since $ \int_x^{\lambda x} \frac{g'(y)}{g(y)}\, d y= \ln
(g(\lambda x)) -\ln(g(x))$, the result is proved.
\end{pf}
\begin{lemma}
\label{lemmaphi} Suppose that $\phi$ is regularly varying at $\infty
$ with index $\beta\in[0,1[$ and that $\phi(x)\rightarrow\infty$
as $x \rightarrow\infty$. Let $\beta' \in\,]\beta,1[$. There then
exists some $x_1(\beta')$ such that for $x \geq x_1(\beta')$
and all $\lambda\geq1$,
\[
1 \leq\phi(x) \leq\phi(\lambda x) \leq\lambda^{\beta'}\phi(x)
\]
and
\[
\phi(x) \leq x^{\beta'}.
\]
\end{lemma}
\begin{pf}
Note that $\phi$ is infinitely
differentiable on $]0,\infty[$ with derivative
\[
\phi'(x)=\int_0^{\infty} v\exp(-x v) \Pi( d v),
\]
which is nonincreasing. It is then a classical result on regular
variation (see the monotone density theorem, \cite{BGT}, Theorem
1.7.2) that $\phi'$ is regularly varying with index $\beta-1$ and
\[
\frac{x\phi'(x)}{\phi(x)}\mathop{\rightarrow}_{x \rightarrow
\infty} \beta.
\]
The first part of the lemma is then a consequence of the above Lemma
\ref{tech1} and of the fact that $\phi$ is increasing and converges
to $\infty$. We also have that $\phi(x)/x^{\beta'}$ converges to $0$
at $\infty$ (since $\beta'>\beta$), hence the second assertion holds
for sufficiently large $x$.
\end{pf}
\begin{lemma}
\label{tech3}
Let
\[
f(x)=-\ln\bigl(\mathbb P(I>x) \bigr),\qquad x\geq0,
\]
which, as proved in Lemma \ref{Victor}, is regularly varying with
index $1/(1-\beta)$, under the assumption (\ref{hypothese}). Suppose,
moreover, that $\int^{\infty}x \Pi( d x)<\infty$. Then $f$
is infinitely differentiable and
\[
\frac{xf'(x)}{f(x)} \rightarrow\frac{1}{1-\beta} \qquad\mbox{as } x
\rightarrow\infty.
\]
\end{lemma}
\begin{pf}
According to \cite{CarmonaPetitYor},
Proposition 2.1, when $\int^{\infty}x \Pi( d
x)<\infty$, there exists an infinitely differentiable function $k \dvtx
]0,\infty[\  \rightarrow[0,\infty[$ such that $k(x)\, dx$ is the
distribution of $I$. Moreover,
\begin{eqnarray*}
k(x)&=&\int_x^{\infty} \overline{\Pi} \bigl(|\alpha|^{-1} \ln
(u/x) \bigr) k(u) \, du \\
&=&\int_0^{\infty} \biggl( \int_{x}^{xe^{v|\alpha|}}k(u) \, d
u \biggr) \Pi
( d v).
\end{eqnarray*}
To simplify notation, we suppose in the following that $|\alpha|=1$.
The proof is identical for $|\alpha|\neq1$.
In particular, we have
\[
\mathbb P(I>x)=\int_x^{\infty} k(u) \, du
\]
and
%
%
\begin{equation}
\label{fprime}\qquad
f'(x)=\frac{k(x)}{\mathbb P(I>x)}=\int_0^{\infty} \bigl(1-\exp
\bigl(f(x)-f(xe^{v})
\bigr) \bigr) \Pi( d v),\qquad x>0.
\end{equation}
Note that since $f$ is regularly varying with a positive index, we have
that $f(x) \rightarrow\infty$ as $x \rightarrow\infty$ and,
therefore, for all $v>0$,
\[
f(x)-f(xe^v)=f(x) \biggl(1-\frac{f(xe^v)}{f(x)} \biggr) \mathop{\sim
}_{\infty} f(x) \bigl(1-e^{v/(1-\beta)} \bigr)\mathop{\rightarrow
}_{x \rightarrow\infty} -\infty.
\]

\begin{itemize}
\item When $\Pi(]0,\infty[)<\infty$, this implies the
expected result since, by dominated convergence,
\[
f'(x) \mathop{\rightarrow}_{x \rightarrow\infty} \Pi(]0,\infty[)
=\lim_{x \rightarrow\infty} \frac{f(x)}{x}.
\]
\item The proof is much more technical when $\Pi
(]0,\infty[)=\infty$, which is supposed for the rest of this proof.
We proceed in two steps.
\end{itemize}

\textit{Step} 1. The goal of this step is to prove that
\[
\liminf_{x \rightarrow\infty} \frac{xf'(x)}{f(x)} \geq\frac
{1}{1-\beta}.
\]
First, suppose that there exists some $x_0$ and some \textit
{nondecreasing} positive function $g$ such that $f'(x) \geq g(x)$ for
all $x \geq x_0$. Then, for $x \geq x_0$ and $v>0$,
\[
f(xe^v)-f(x) = \int_x^{xe^v} f'(u) \, du \geq
g(x)x(e^v-1) \geq g(x)x v.
\]
Using (\ref{fprime}), this gives
%
%
\begin{equation}
\label{iterate}
f'(x) \geq\phi(g(x)x),\qquad x \geq x_0.
\end{equation}
Now, note that $f'(x) \rightarrow\infty$ as $x \rightarrow\infty$
since, for all $a>0$,
\[
\liminf_{x \rightarrow\infty} f'(x) \geq
\liminf_{x \rightarrow\infty} \int_a^{\infty}
\bigl(1-\exp\bigl(f(x)-f(xe^v) \bigr) \bigr)\Pi( d v)
= \int_a^{\infty}\Pi( d v)
\]
(by dominated convergence) and the right-hand side converges to $\infty
$ as $a\rightarrow0$. In particular, $f'(x) \geq1$ for $x$
sufficiently large (say $x \geq x_0$).
Replacing $g$ by $1$ in (\ref{iterate}), we get
\[
f'(x) \geq\phi(x)\qquad \forall x \geq x_0.
\]
Recall that $\phi$ is nondecreasing and then iterate the procedure to
get, for all $n \geq0$,
%
%
\begin{equation}
\label{ineg5}
f'(x) \geq h_n(x)\qquad \forall x \geq x_0,
\end{equation}
where the functions $h_n \dvtx]0,\infty[\  \rightarrow\mathbb]0,\infty[$
are defined by induction by
\begin{eqnarray*}
&& h_0(x)=1\qquad \mbox{for all } x\geq0; \\
&& h_n(x)=\phi(h_{n-1}(x) x)\qquad \mbox{for all } x\geq0.
\end{eqnarray*}
Now, the interesting fact is that for $x$ large enough, $h_n(x)
\rightarrow\varphi(x)/x$ as $n \rightarrow\infty$. Indeed, let
$\beta' \in\,]\beta,1[$. With the notation of Lemma \ref{lemmaphi},
we have, for $x \geq x_1(\beta')$, $1 \leq\phi(x) \leq x^{\beta'}$,
that is, $ h_0(x)\leq h_1(x) \leq x^{\beta'}$. Using the fact that
$\phi$ is nondecreasing, we easily have, by induction, that
\[
1 \leq h_n(x) \leq h_{n+1}(x) \leq x^{\beta'+\cdots+\beta'^{n+1}}
\leq x^{\beta'/(1-\beta')}<\infty
\]
for all $n \geq1$. Let $l(x):=\lim_{n \rightarrow\infty} h_n(x)$.
We have shown that $0<l(x)<\infty$.
Then, necessarily, $l(x)=\phi(l(x)x)$ [in other words, $l(x)x/\phi
(l(x)x)=x$] and, finally, $l(x)x=\varphi(x)$, $\forall x \geq
x_1(\beta')$.
To conclude, for $x$ large enough, letting $n \rightarrow\infty$
in~(\ref{ineg5}), we get
$
f'(x) \geq\varphi(x)/x,
$
which, combined with Lemma \ref{Victor}, gives the expected lim inf.

\textit{Step} 2. The proof of the lim sup is similar, but
more technical. First, note that for all $\varepsilon>0$ and all
$a<\ln(1+\varepsilon)$, $a>0$,
\begin{eqnarray*}
&&\liminf_{x \rightarrow\infty} \int_0^{\ln(1+\varepsilon)}
\bigl(1-\exp\bigl(f(x)-f(xe^v) \bigr) \bigr)\Pi( d v)
\\
&&\qquad\geq \liminf_{x \rightarrow\infty} \int_a^{\ln(1+\varepsilon)}
\bigl(1-\exp\bigl(f(x)-f(xe^v) \bigr) \bigr)\Pi( d v) \\
&&\qquad= \int_a^{\ln(1+\varepsilon)}\Pi( d v)\qquad \mbox{(by
dominated convergence),}\\
&&\qquad\hspace*{-2.91pt}\mathop{\rightarrow}_{a \rightarrow0} \infty,
\end{eqnarray*}
whereas
\[
\lim_{n\rightarrow\infty}\int_{\ln(1+\varepsilon)}^{\infty}
\bigl(1-\exp\bigl(f(x)-f(xe^v) \bigr) \bigr)\Pi( d
v)=\int_{\ln(1+\varepsilon)}^{\infty}\Pi( d
v)<\infty.
\]
Hence, there exists some $x_1(\varepsilon)$ such that for $x \geq
x_1(\varepsilon)$,
%
%
\begin{equation}
\label{inegfprime}
f'(x)\leq(1+\varepsilon)\int_0^{\ln(1+\varepsilon)}
\bigl(1-\exp\bigl(f(x)-f(xe^v) \bigr) \bigr)\Pi( d v).
\end{equation}
Next, fix some $\beta' \in\,]\beta,1[$ and consider some $\delta>0$
and $\varepsilon>0$ such that
$(1+\delta)(1+\varepsilon)^{1/(\beta-1)}\beta'<1$. Since $f$ is
regularly varying with index $1/(1-\beta)$, there exists some
$x_2(\delta,\varepsilon)$ such that
%
%
\begin{equation}
\label{inegf}
f\bigl(x(1+\varepsilon)\bigr)\leq(1+\delta)(1+\varepsilon)^{1/(1-\beta
)}f(x)\qquad \forall x \geq x_2(\delta,\varepsilon).
\end{equation}
We will need this later. For the moment, let $x_0=\max(x_1(\beta
'),x_1(\varepsilon), x_2(\delta,\varepsilon))$, with $x_1(\beta')$
as introduced in Lemma \ref{lemmaphi}.
Next, suppose that for all $x \geq x_0$,
\[
f'(x) \leq g(x)
\]
for some \textit{nondecreasing} function $g$ such that $g(x) \geq1$
for all $x \geq x_0$. Note that this implies that
\[
f(xe^v)-f(x) = \int_x^{xe^v} f'(u) \, du \leq
g(xe^v)x(e^v-1).
\]
The function $v \rightarrow v^{-1}(e^v-1)$ is
increasing on $[0,\infty[$, hence $e^v-1 \leq v \gamma(\varepsilon)$
for all $v
\leq\ln(1+\varepsilon)$, where
$\gamma(\varepsilon)=\varepsilon/(\ln(1+\varepsilon))$.
Together with (\ref{inegfprime}),
this leads to
%
%
\begin{eqnarray}
\label{inegfprime2}
f'(x)&\leq&(1+\varepsilon)
\int_0^{\ln(1+\varepsilon)} \bigl(1-\exp\bigl(-
g\bigl(x(1+\varepsilon)\bigr)x v\gamma(\varepsilon) \bigr) \bigr) \Pi
( d v)\nonumber\\[-8pt]\\[-8pt]
&\leq&(1+\varepsilon) \phi\bigl( g\bigl(x(1+\varepsilon)\bigr)x
\gamma(\varepsilon) \bigr)\nonumber
\end{eqnarray}
for all $x \geq x_0$.

We then claim that for all $n \geq1$ and all $x \geq x_0$,
%
%
\begin{eqnarray}\hspace*{26pt}
\label{recfprime}
f'(x) \leq(1+\varepsilon)^{1+\beta'+2\beta'^2+\cdots+n \beta'^n}
\gamma(\varepsilon)^{\beta'+\beta'^2+\cdots+\beta'^n}g
\bigl(x(1+\varepsilon)^n \bigr)^{\beta'^n} h_n(x),
\end{eqnarray}
where the sequence of functions $h_n$ is that introduced in step 1 of
this proof.
We will prove this by induction on $n$. First, though, let us mention
that, by a simple application of induction, using Lemma \ref{lemmaphi},
\[
h_n\bigl(x(1+\varepsilon)\bigr) \leq(1+\varepsilon)^{\beta'+\cdots+\beta
'^n} h_n(x)\qquad\mbox{for all } x \geq x_0 \mbox{ and } n \geq1.
\]
We now turn to the proof of (\ref{recfprime}). For $n=1$, we can use
(\ref{inegfprime2}) and Lemma \ref{lemmaphi} to get [note that
$\gamma(\varepsilon) \geq1$, hence $\gamma(\varepsilon)g \geq1$]
\[
f'(x) \leq(1+\varepsilon) g\bigl(x(1+\varepsilon)\bigr)^{\beta'}\gamma
(\varepsilon)^{\beta'} \phi(x),\qquad x \geq x_0,
\]
which leads to (\ref{recfprime}) for $n=1$. Now, assume that (\ref
{recfprime}) is true for some integer $n$. Note that the function on
the right-hand side of this inequality, which we call $g_1$, is larger
than $1$ for all $x \geq x_0$. Also, note that it is nondecreasing.
Hence, we get, replacing $g$ by $g_1$ in (\ref{inegfprime2}), for $x
\geq x_0$,
\begin{eqnarray*}
f'(x) &\leq&(1+\varepsilon) \phi\bigl( (1+\varepsilon)^{1+\beta
'+2\beta'^2+\cdots+n \beta'^n} \gamma(\varepsilon)^{\beta'+\beta
'^2+\cdots+\beta'^n}\\
&&\hspace*{50.26pt}{}\times g \bigl(x(1+\varepsilon)^{n+1} \bigr)^{\beta'^n}
h_n\bigl(x(1+\varepsilon)\bigr) x \gamma(\varepsilon) \bigr) \\
&\leq& (1+\varepsilon) \phi\bigl( (1+\varepsilon)^{1+2\beta
'+3\beta'^2+\cdots+(n+1) \beta'^n}\\
&&\hspace*{40.71pt}{}\times  \gamma(\varepsilon)^{1+\beta
'+\beta'^2+\cdots+\beta'^n}g \bigl(x(1+\varepsilon)^{n+1}
\bigr)^{\beta'^n} h_n(x) x \bigr) \\
&\leq& (1+\varepsilon)^{1+\beta'+2\beta'^2+\cdots+(n+1) \beta'^{n+1}}
\gamma(\varepsilon)^{\beta'+\beta'^2+\cdots+\beta'^{n+1}}\\
&&{}\times
g \bigl(x(1+\varepsilon)^{n+1} \bigr)^{\beta'^{n+1}} \phi(
h_n(x) x ),
\end{eqnarray*}
where, for the last inequality, we have used Lemma \ref{lemmaphi}.
Hence, we have (\ref{recfprime}) for all $n \geq1$.

Now, thanks to the assumptions (\ref{hypothese}) and $\int^{\infty}x
\Pi( dx)<\infty$ and to Lemma 1 of~\cite{HaasRegularity}, we
know that the function $k$ is bounded from above on $]0,\infty[$, say
by some constant $C \geq1$. Hence, $f'(x)=k(x)/\mathbb P(I>x) \leq C
\exp(f(x))$ for all \mbox{$x > 0$}. Since $f$ is nondecreasing and
nonnegative, the function $x \rightarrow C \exp(f(x))$ is
nondecreasing and greater than $1$, hence we can replace $g$ by this
function in (\ref{recfprime}) to get, for all $n \geq1$ and all $x
\geq x_0$,
%
%
\begin{eqnarray}
\label{ineg3fprime}
f'(x) &\leq& (1+\varepsilon)^{1+\beta'+2\beta'^2+\cdots+n \beta'^n}
\gamma(\varepsilon)^{\beta'+\beta'^2+\cdots+\beta'^n}C^{\beta
'^n}\nonumber\\[-8pt]\\[-8pt]
&&{}\times \exp\bigl(\beta'^n f\bigl(x(1+\varepsilon)^n\bigr) \bigr) h_n(x).\nonumber
\end{eqnarray}
Our goal now is to let $n \rightarrow\infty$ in this inequality.
Iterating inequality (\ref{inegf}), we get, for $x \geq x_0$ and for
all $n \geq1$,
\[
f\bigl(x(1+\varepsilon)^n\bigr)\leq(1+\delta)^n(1+\varepsilon)^{n/(1-\beta)}f(x).
\]
Since
\[
(1+\delta)(1+\varepsilon)^{1/(\beta-1)}\beta'<1,
\]
this leads, for $x \geq x_0$, to
\[
\exp\bigl(\beta'^n f\bigl(x(1+\varepsilon)^n\bigr) \bigr) \mathop
{\rightarrow}_{n \rightarrow\infty} 1.
\]
As $n \rightarrow\infty$, we also have
\[
C^{\beta'^n} \rightarrow1 \quad\mbox{and}\quad (1+\varepsilon)^{1+\beta
'+2\beta'^2+\cdots+n \beta'^n} \rightarrow(1+\varepsilon)^{1+\beta
'/(1-\beta')^2}
\]
and
\[
\gamma(\varepsilon)^{\beta'+\beta
'^2+\cdots+\beta'^n} \rightarrow\gamma(\varepsilon)^{\beta'/(1-\beta')}.
\]
Last, recall that for $x$ large enough, $h_n(x) \rightarrow\varphi
(x)/x$ as $n \rightarrow\infty$.
Letting $n \rightarrow\infty$ in (\ref{ineg3fprime}), we therefore
have, for $x$ large enough,
\[
f'(x) \leq C_{\varepsilon} \varphi(x)/x,
\]
where $C_{\varepsilon} \rightarrow1$ as $\varepsilon\rightarrow0$.
This gives
\[
\limsup_{x \rightarrow\infty} \frac{xf'(x)}{f(x)} \leq\frac
{1}{1-\beta}.
\]
\upqed\end{pf}
\begin{lemma}
\label{tech4}
Let
\[
f(x):=-\ln(m(x) ),\qquad x\geq0,
\]
and suppose that (\ref{hypothese}) holds, $\int^{\infty}x \Pi
( d x)<\infty$ and $\mu_0$ has bounded support. Then $f$ is
differentiable on $]0,\infty[$ and
\[
\frac{xf'(x)}{f(x)}=-\frac{xm'(x)}{m(x)f(x)} \rightarrow\frac
{1}{1-\beta} \qquad\mbox{as } x \rightarrow\infty.
\]
\end{lemma}
\begin{pf}
With no loss of generality, we suppose that the supremum of the support
of $\mu_0$ is equal to $1$. Under the assumptions of the lemma, we
know (see the proof of the previous lemma) that $x \rightarrow\mathbb
P(I>x) $ is differentiable on $]0,\infty[$, with derivative $-k$. By
Lemma 1 in \cite{HaasRegularity}, we also know that the function $x
\in\,]0,\infty[\  \rightarrow xk(x)$ is bounded. Let $M$ denote an upper
bound. Recall, then, that
\[
m(x)=\int_0^1 \mathbb P(I>xy^{\alpha}) y \mu_0( d y)
\]
and note that for all $x>a>0$ and all $y \in\,]0,1[$,
\[
\bigl|\partial_x \bigl( \mathbb P(I>xy^{\alpha}) \bigr)
\bigr|=k(xy^{\alpha})y^{\alpha} \leq\frac{M}{a}.
\]
Hence, by dominated convergence, $m$ is continuously differentiable on
$]0,\infty[$, with derivative
\[
m'(x)=-\int_0^1 k(xy^{\alpha}) y^{\alpha}y \mu_0( d y),\qquad x>0.
\]
Now, fix $\delta>0$. By Lemma \ref{tech3}, there exists some
$x(\delta)$ such that for $x \geq x(\delta)$,
\[
\frac{1-\delta}{1-\beta} \leq\frac{-xk(x)}{\mathbb P(I>x) \ln
(\mathbb P(I>x) )} \leq\frac{1+\delta}{1-\beta}.
\]
Then, for $x \geq x(\delta)$,
%
%
\begin{eqnarray}
\label{inter1}
&&\frac{1-\delta}{(1-\beta)x} \int_0^{1} \mathbb P(I>xy^{\alpha})
\ln\bigl( \mathbb P(I>xy^{\alpha}) \bigr) y \mu_0( d y)
\nonumber\\[-8pt]\\[-8pt]
&&\qquad
\leq m'(x) \leq\frac{1+\delta}{(1-\beta)x} \int_0^{1} \mathbb
P(I>xy^{\alpha}) \ln\bigl( \mathbb P(I>xy^{\alpha}) \bigr) y \mu
_0( d y).\nonumber
\end{eqnarray}
Now, let $\varepsilon>0$. On the one hand, we claim that
%
%
\begin{eqnarray}
\label{inter2}
&&\int_{1-\varepsilon}^{1} \mathbb P(I>xy^{\alpha}) \ln\bigl(
\mathbb P(I>xy^{\alpha}) \bigr) y \mu_0( d y) \nonumber\\[-8pt]\\[-8pt]
&&\qquad\mathop{\sim
}_{x \rightarrow\infty} \int_{0}^{1} \mathbb P(I>xy^{\alpha}) \ln
\bigl( \mathbb P(I>xy^{\alpha}) \bigr) y \mu_0( d y).\nonumber
\end{eqnarray}
Indeed, for all $0<y<1-\varepsilon$,
\[
\frac{ \mathbb P(I>xy^{\alpha}) \ln( \mathbb P(I>xy^{\alpha
}) )}{ \mathbb P(I>x(1-\varepsilon)^{\alpha}) \ln(
\mathbb P(I>x(1-\varepsilon)^{\alpha}) )} \rightarrow0
\qquad\mbox{as } x\rightarrow\infty
\]
since $x \rightarrow-\ln( \mathbb P(I>x) )$ is regularly
varying at $\infty$ with a positive index and $\alpha<0$. It is then
not hard to see, using Lemmas \ref{tech1} and \ref{tech3}, that for
$x$ large enough, this function is bounded from above by
\[
\exp\biggl( 1- \biggl(\frac{y}{1-\varepsilon} \biggr)^{{\alpha
(1-\varepsilon)}/({1-\beta})} \biggr) \biggl(\frac{y}{1-\varepsilon}
\biggr)^{{\alpha(1+\varepsilon)}/({1-\beta})},
\]
which, in turn, is bounded for $y \in\,]0,1-\varepsilon[$. Hence, by
dominated convergence, we see that
\begin{eqnarray*}
&&\biggl| \int_{0}^{1-\varepsilon} \mathbb P(I>xy^{\alpha}) \ln
\bigl( \mathbb P(I>xy^{\alpha}) \bigr) y \mu_0( d y)
\biggr|\\
&&\qquad\mathop{\ll}_{x \rightarrow\infty}
\bigl|\mathbb P\bigl(I>x(1-\varepsilon)^{\alpha}\bigr) \ln\bigl( \mathbb
P\bigl(I>x(1-\varepsilon)^{\alpha}\bigr) \bigr) \bigr| \\
&&\hspace*{5.22pt}\qquad\leq \frac
{1}{\varepsilon} \biggl| \int_{1-\varepsilon}^{1} \mathbb
P(I>xy^{\alpha}) \ln\bigl( \mathbb P(I>xy^{\alpha}) \bigr) y \mu
_0( d y) \biggr|,
\end{eqnarray*}
where, for the last inequality, we have used the fact that the function
$x \rightarrow-x\ln(x)$ is increasing in a neighborhood of $0$. Hence
(\ref{inter2}).
A similar, but simpler, argument leads to the result
%
%
\begin{equation}
\label{inter3}
\int_{1-\varepsilon}^{1} \mathbb P(I>xy^{\alpha}) y \mu_0( d
y) \mathop{\sim}_{x \rightarrow\infty} \int_{0}^{1} \mathbb
P(I>xy^{\alpha}) y \mu_0( d y).
\end{equation}
On the other hand, using the fact that $x \rightarrow\ln
(\mathbb P(I>x) )$ is regularly varying with index $1/(1-\beta
)$, we have, for $1-\varepsilon\leq y \leq1$ and $x$ sufficiently
large [say $x \geq x(\varepsilon)$],
\begin{eqnarray*}
(1+\varepsilon)(1-\varepsilon)^{\alpha/(1-\beta)}\ln\bigl(
\mathbb P(I>x) \bigr)
&\leq&\ln\bigl( \mathbb P\bigl(I>x(1-\varepsilon)^{\alpha}\bigr) \bigr) \leq
\ln\bigl( \mathbb P(I>xy^{\alpha}) \bigr)\\
&\leq& \ln\bigl( \mathbb P(I>x) \bigr).
\end{eqnarray*}
Thus,
\begin{eqnarray*}
\int_{1-\varepsilon}^{1} \mathbb P(I>xy^{\alpha}) y \mu_0( d
y) &\leq& \frac{\int_{1-\varepsilon}^{1} \mathbb P(I>xy^{\alpha}) \ln
( \mathbb P(I>xy^{\alpha}) ) y \mu_0( d y)}{\ln
( \mathbb P(I>x) )}\\
&\leq&
(1+\varepsilon)(1-\varepsilon
)^{\alpha/(1-\beta)} \int_{1-\varepsilon}^{1} \mathbb
P(I>xy^{\alpha}) y \mu_0( d y),
\end{eqnarray*}
which, taking $x(\varepsilon)$ larger if necessary and using (\ref
{inter2}) and (\ref{inter3}), gives, for $x \geq x(\varepsilon)$,
\begin{eqnarray*}
(1-\delta) m(x)
&\leq&
\frac{\int_{0}^{1} \mathbb P(I>xy^{\alpha})
\ln( \mathbb P(I>xy^{\alpha}) ) y \mu_0( d
y)}{\ln( \mathbb P(I>x) )}\\
&\leq&
(1+\delta)(1+\varepsilon
)(1-\varepsilon)^{\alpha/(1-\beta)} m(x).
\end{eqnarray*}
Plugging this into (\ref{inter1}) and letting first $\varepsilon
\rightarrow0$ and then $\delta\rightarrow0$, we get the expected
convergence since $f(x) \sim-\ln(\mathbb P(I>x))$ as $x \rightarrow
\infty$.
\end{pf}

\subsubsection{Proof of Theorem \protect\ref{ThCVSub}\textup{(ii)}}
The fact that the function
\[
x \in\,]0,\infty[\  \rightarrow f(x):=-\ln(m(x)) =-\ln\bigl( \mathbb
P\bigl(X(x)>0\bigr) \bigr)
\]
is continuously differentiable on $]0,\infty[$ can be proven in
exactly the same way as when the support of $\mu_0$ is compact; see
the beginning of the proof of Lemma \ref{tech4}.
Next, by Karamata's theorem (Theorem 1.5.11 of \cite{BGT}), if $f$
varies regularly at $\infty$ with index $\lambda>0$ and if its
derivative is also regularly varying at $\infty$, then
\[
\frac{xf'(x)}{f(x)} \rightarrow\lambda\qquad\mbox{as } x \rightarrow
\infty.
\]
Together with Theorem \ref{notcompact}(i), Lemma \ref{tech1} and
Lemma \ref{lemmaCVSub}, this implies Theorem~\ref{ThCVSub}(ii).

\subsection{Quasi-stationary distributions}
\label{QS}

\mbox{}

\begin{pf*}{Proof of Theorem \protect\ref{ThQS}}
When $X(0) \sim\mu_R^{(\lambda)}$,
the distribution of $X(0)^{|\alpha|}I$ is that of $\lambda^{|\alpha
|}R I$, with $R$ independent of $I$,
that is, that of an exponential random variable with parameter
$\lambda^{\alpha}$. We then immediately have that for $n \geq1$ and
$t \geq0$,
\[
\mathbb E\bigl [ \bigl(\bigl(X(0)^{|\alpha|}I-t\bigr)^+ \bigr)^n
\bigr]=\lambda^{|\alpha|
n}n!\exp(-\lambda^{\alpha} t)
\]
and
\[
\mathbb P\bigl(X(t)>0\bigr)=\mathbb P\bigl(X(0)^{|\alpha|}I>t\bigr)=\exp(-\lambda
^{\alpha}t).
\]
Following the beginning of the proof of Lemma \ref{lemmaCVSub}, this gives
\[
\mathbb E \bigl[ ( X(t)
)^{|\alpha| n} \bigr]\mathbb E[I^n] =\mathbb E \bigl[
\bigl(\bigl(X(0)^{|\alpha|}I-t\bigr)^+ \bigr)^n \bigr]
=\lambda^{|\alpha| n}n!\exp(-\lambda^{\alpha} t)
\]
and then
\[
\mathbb E \bigl[ ( X(t)
)^{|\alpha| n}| X(t)>0 \bigr]=\mathbb E\bigl[\lambda^{|\alpha|
n}R^n\bigr] =\mathbb E \bigl[ X(0)^{|\alpha| n} \bigr].
\]
Hence, $\mu_R^{(\lambda)}$ is a quasi-stationary distribution since
the distribution of
$R$ is characterized by its entire positive moments. Note that there is
no other quasi-stationary distribution. Indeed, let $\varsigma$ be a
quasi-stationary distribution and suppose that $X(0) \sim\varsigma$.
Then, necessarily, by the Markov property of $X$, $\mathbb
P(X(t+s)>0)=\mathbb P(X(t)>0)\mathbb P(X(s)>0)$, which implies that
$X(0)^{|\alpha|}I$ has an exponential distribution, say with parameter
$\ell$, that is, $\ell X(0)^{|\alpha|}I$ has an exponential
distribution with parameter $1$. Since the factorization (\ref{face})
characterizes the distribution of $R$, we get that $\varsigma=\mu
_R^{(\ell^{1/\alpha})}$.
\end{pf*}
\begin{pf*}{Proof of Theorem \protect\ref{stat}}
The first part of this theorem is an obvious consequence of Theorem \ref{ThQS}. The
reverse cannot be directly deduced from Theorem \ref{ThQS} since we do
not know if uniqueness holds for the fragmentation equation when the
initial measure has an unbounded support.

So, consider $(\mu_t,t \geq0)$, a quasi-stationary solution of the
fragmentation equation (\ref{eqfrag}). We want to prove that this
solution belongs to the family of solutions $ ((\mu^{(\lambda
)}_{\infty,t}, t \geq0), \lambda> 0 )$, as defined in Theorem
\ref{stat}. Replacing $\mu_t$ by $m(t)\mu_0$ in equation (\ref
{eqfrag}), we get that
\[
\bigl(1-m(t)\bigr) \langle\mu_0,f\rangle=-\int_0^t m(s) \, ds \langle\mu_0,
G(f)\rangle\qquad \forall
f \in C^1_c,
\]
where $G(f)(x)=x^{\alpha}\int_0^1 (f(xy)-f(x)y)B( dy)$.
In other words, there exists some constant $C>0$ such that
\[
m(t)=\exp(-Ct)\qquad \forall t \geq0,
\]
and
\[
\langle\mu_0,f\rangle=-C^{-1}\langle\mu_0, G(f)\rangle\qquad \forall f \in C^1_c.
\]
When $f \in C^1_c$, the function $x\rightarrow xf(x)$ is also in
$C^1_c$. Hence, the above identity can be rewritten
%
%
\begin{equation}
\label{C}
\langle x\mu_0,f\rangle=-C^{-1}\langle x\mu_0, \tilde A(f)\rangle\qquad
\forall f \in C^1_c,
\end{equation}
where $\tilde A(f)(x)=x^{\alpha}\int_0^1 (f(xy)-f(x))yB( dy)$.

To show that this characterizes $\mu_0$, we need the following fact:
for all $\beta>0$, there exists a nondecreasing sequence of functions
$f_{\beta,n}\dvtx]0,\infty[\  \rightarrow[0,\infty[$ such that $f_{\beta
,n}(x)\rightarrow x^{\beta}$ as $n \rightarrow\infty$, $\forall
x>0$, $f_{\beta,n} \in C^1_c$ and $|f'_{\beta,n}(x)| \leq\beta
x^{\beta-1}$ for all $x>0$ and all $n \geq1$. This sequence can, for
example, be constructed by first considering a nondecreasing sequence
of continuous functions $g_{\beta,n}\dvtx]0,\infty[\  \rightarrow[0,\infty
[$ such that $g_{\beta,n}(x) \leq\beta x^{\beta-1}$, $\forall x>0, n
\geq1$, $g_{\beta,n}(x)=\beta x^{\beta-1}$ for $x \in[n^{-1},n]$
and $g_{\beta,n}(x)=0$ for $x \in\,]0,(2n)^{-1}] \cup[2n,\infty[$.
Then, set $f_{\beta,n}(x):=\int_0^{x} g_{\beta,n}(u) \, du$ for
$x \in\,]0,2n]$ and extend these functions to $]2n,\infty[$ so that
$f_{\beta,n} \in C^1_c$ and $|f'_{\beta,n}(x)| \leq\beta x^{\beta
-1}$, for all $x>0$ and all $n \geq1$, and the sequence $(f_{\beta
,n},n\geq1)$ is nondecreasing.
For all $\beta>0$, this implies that for all $x>0$,
%
%
\begin{equation}
\label{A}
\tilde A(f_{\beta,n})(x) \mathop{\rightarrow}_{n \rightarrow\infty
} x^{\alpha+\beta} \int_0^1 (y^{\beta}-1)y B( dy)=-x^{\alpha
+\beta} \phi(\beta),
\end{equation}
together with
%
%
\begin{equation}
\label{B}
|\tilde A(f_{\beta,n})(x) | \leq(2+\beta) x^{\alpha+\beta} \int
_0^1(1-y)y B( dy).
\end{equation}
Indeed, this is obvious when $\beta\geq1$: we just need that
\begin{eqnarray*}
|f_{\beta,n}(xy)-f_{\beta,n}(x)| &\leq& {\sup_{z \in[xy,x]}}|f'_{\beta
,n}(z)| x(1-y) \\
&\leq& \beta x^{\beta}(1-y)\qquad \mbox{for } y \in\,]0,1[, x>0,
\end{eqnarray*}
and then use the dominated convergence theorem. The case $0<\beta<1$
needs more care. Using the aforementioned properties of $f_{\beta,n}$
and also the fact that $f_{\beta,n}(x) \leq x^{\beta}$, we obtain,
for $x>0$ and $y \in\,]0,1[$,
\begin{eqnarray*}
&&
x|f_{\beta,n}(xy)- f_{\beta,n}(x)|\\
&&\qquad\leq x f_{\beta,n}(xy) (1-y)+
|xyf_{\beta,n}(xy)- xf_{\beta,n}(x)| \\
&&\qquad \leq x^{1+\beta}(1-y) + {\sup_{z \in[xy,x]}}|(\mathrm{id}
f_{\beta,n})'(z)| x(1-y) \\
&&\qquad\leq x^{1+\beta}(1-y) +(1+\beta) x^{1+\beta}(1-y),
\end{eqnarray*}
which leads to (\ref{A}) and (\ref{B}).

Now, take $\beta=|\alpha|$. Then use (\ref{A}), (\ref{B}) and the
dominated convergence theorem on the right-hand side of (\ref{C})
[recall that $x\mu_0( dx)$ is a probability measure], together
with the monotone convergence theorem on the left-hand side of (\ref
{C}), to get
\[
\int_0^{\infty} x^{|\alpha|} x\mu_0( dx)=C^{-1} \phi
(|\alpha|)<\infty.
\]
Then, by an obvious induction, taking successively $\beta=2|\alpha|$,
$\beta=3|\alpha|$, etc., we get, for all $n \geq1$, that
\begin{eqnarray*}
\int_0^{\infty} x^{n|\alpha|} x\mu_0( dx)
&=& C^{-1} \phi
(n|\alpha|)\int_0^{\infty} x^{(n-1)|\alpha|} x\mu_0(
dx)\\
&=& C^{-n} \phi(n|\alpha|)\cdots\phi(|\alpha|).
\end{eqnarray*}
We recognize the moments formula (\ref{defmuR}). Hence, $x \mu
_0( dx)=\mu_R^{(C^{1/\alpha})}( d x)=x\mu_{\infty
}^{(C^{1/\alpha})}( dx)$ and for all $t \geq0$, $\mu
_t=m(t)\mu_0=\exp(-Ct)\mu_{\infty}^{(C^{1/\alpha})}=\mu_{\infty
,t}^{(C^{1/\alpha})}$.
\end{pf*}

\section{Different speeds of decrease: proof of Proposition \protect
\ref{fasterslower}}
\label{SectionFS}

\subsection{Proof of Proposition
\protect\ref{fasterslower}\textup{(i)}}

Recall that the support of $\mu_0$ is supposed to be bounded with
supremum $1$. The goal of this section is to prove the forthcoming
Corollary \ref{corofasterslower}, which is the statement of
Proposition \ref{fasterslower}(i) translated in terms of the process
$X$ defined by (\ref{defX}), provided the L\'evy measure $\Pi$ of the
subordinator $\xi$ involved in the construction of $X$ is related to
the fragmentation measure $B$ by (\ref{PiB}) and $X(0)$ is distributed
according to $x\mu_0( dx)$. We recall that the distribution of
$X(t)$ conditional on $X(t)>0$ is then $x\mu_t( dx)/m(t)$, $t
\geq0$, where $(\mu_t,t\geq0)$ denotes the solution of the
fragmentation equation starting from $\mu_0$.
We start with some preliminary lemmas.
\begin{lemma}
\label{Lemma2} Suppose that (\ref{hypothese}) holds and that $\int_0
|{\ln(x)}|x B( dx)<\infty$. Consider some random variable $I$
independent of $X$, with distribution that of $\int_0^{\infty} \exp
(\alpha\xi_r)\, d r$. Then:

\begin{longlist}
\item there exists some $t_0>0$ such that
\[
\sup_{t \geq t_0, a>0} a^{\alpha}\mathbb P \biggl( \biggl( \frac
{\varphi(|\alpha|t)}{|\alpha|t} \biggr)^{1/|\alpha|}X(t)
I^{1/|\alpha|} \leq a \Big| X(t)>0 \biggr) <\infty;
\]
\item for all positive functions $g\dvtx[0,\infty[\  \rightarrow\,]0,\infty[
$ converging to $0$ at $\infty$, we have, as $t \rightarrow\infty$,
\[
g(t)^{\alpha}\mathbb P \biggl( \biggl( \frac{\varphi(|\alpha
|t)}{|\alpha|t} \biggr)^{1/|\alpha|}X(t) I^{1/|\alpha|} \leq g(t)
\Big| X(t)>0 \biggr) \rightarrow1.
\]
\end{longlist}
\end{lemma}
\begin{pf}
To simplify notation, suppose that $\alpha
=-1$ (the proof is identical for all $\alpha<0$). Recall, then, the
key equality in law (\ref{keyidentity}), which leads to the following
identities for all $a>0$:
%
%
\begin{eqnarray}
\label{identity1}
&&
a^{-1}\mathbb P \biggl( \frac{\varphi(t)}{t} X(t) I \leq a
| X(t)>0 \biggr)\nonumber\\
&&\qquad= \frac{1}{am(t)}\mathbb P \biggl(
0<X(0)I-t \leq\frac{at}{\varphi(t)} \biggr) \nonumber\\[-8pt]\\[-8pt]
&&\qquad= \frac{m(t)-m(t+at/\varphi(t))}{am(t)} \nonumber\\
&&\qquad= \frac{1}{a} \biggl(1- \exp\biggl( -\ln(m(t)) \biggl( 1-
\frac{\ln(m(t(1+a/\varphi(t))) )}{\ln(m(t))} \biggr)\biggr)
\biggr).\nonumber
\end{eqnarray}
We then use the regular variation of $-\ln( m )$ with
index $1/(1-\beta)$ (Proposition~\ref{mproperties}) and Lemma \ref
{lemmanormalized} to see that for all $\varepsilon>0$, there exists a
real number $t(\varepsilon)$ such that for all $t\geq t(\varepsilon)$
and all $a>0$,
%
%
\begin{eqnarray}
\label{inequalities2}
1- \bigl( 1+a/\varphi(t) \bigr)^{{1}/({1-\beta})+\varepsilon}
&\leq& 1- \frac{\ln(m(t(1+a/\varphi(t))) )}{\ln(m(t))}
\nonumber\\[-8pt]\\[-8pt]
&\leq& 1- \bigl( 1+a/\varphi(t) \bigr)^{{1}/({1-\beta
})-\varepsilon}.\nonumber
\end{eqnarray}
Now, let $0<\varepsilon<1-\beta$. Since
\[
1-(1+x)^{{1}/({1-\beta})+\varepsilon} \geq-x \biggl(\frac
{1}{1-\beta}+2\varepsilon\biggr)
\]
for all $x>0$ sufficiently small and since, further, $\varphi(t)
\rightarrow\infty$ as $t \rightarrow\infty$ and $-\ln(m) \mathop
{\sim}_{\infty} (1-\beta) \varphi$, we have that for all $0<a \leq
1$ and all $t \geq t'(\varepsilon)$ [for some $t'(\varepsilon)$
depending on $\varepsilon$ but not on $0<a \leq1$],
\begin{eqnarray*}
&&
-\ln(m(t)) \bigl(1- \bigl( 1+a/\varphi(t) \bigr)^{
{1}/({1-\beta})+\varepsilon} \bigr)\\
&&\qquad\geq \frac{\ln(m(t))}{\varphi
(t)} a \biggl(\frac{1}{1-\beta}+2\varepsilon\biggr) \\
&&\qquad\geq -a(1-\beta+\varepsilon) \biggl(\frac{1}{1-\beta
}+2\varepsilon\biggr).
\end{eqnarray*}
Together with identities (\ref{identity1}) and inequalities (\ref
{inequalities2}), this implies that for all $t \geq\max(t(\varepsilon
),t'(\varepsilon))$,
\[
\sup_{0< a \leq1}a^{-1}\mathbb P \biggl( \frac{\varphi(t)}{t} X(t)
I \leq a \Big| X(t)>0 \biggr) < \infty.
\]
This is enough to get (i) since for $a \geq1$, $a^{-1}$ multiplied by
a probability is bounded by $1$.

The proof of (ii) relies on the same idea. Since $g(t)/\varphi(t)
\rightarrow0$ as $t \rightarrow\infty$,
\[
-\ln(m(t)) \bigl(1- \bigl( 1+ g(t)/\varphi(t) \bigr)^{
{1}/({1-\beta})+\varepsilon} \bigr) \mathop{\sim}_{\infty}
-g(t)(1-\beta) \biggl(\frac{1}{1-\beta}+\varepsilon\biggr)
\]
and a similar result holds by replacing $\varepsilon$ by $-\varepsilon
$. Together with the inequalities (\ref{inequalities2}) and the
identities (\ref{identity1}) [there replacing $a$ by $g(t)$], also
using the fact that $g(t) \rightarrow0$ as $t \rightarrow\infty$,
we get (ii).
\end{pf}
\begin{lemma}
\label{lemmak}
Suppose that $\kappa:=\int_0^1 |{\ln(x)}| x B( dx)<\infty$. Then:
\begin{longlist}
\item $I$ possesses a density $k \in\mathcaligr C^{\infty
}(]0,\infty[)$;
\item $\mathbb E[I^{-1}]=\kappa|\alpha|<\infty$;
\item if, further, the support of $B$ is not included in a set
of the form $\{ a^n, n \in\mathbb N\}$ for some $a \in\,]0,1[$, then
the function
\[
x \in\mathbb R \rightarrow\mathbb E[I^{ix-1}]=\int_0^{\infty}
y^{ix-1}k(y) \, d y
\]
is well defined and nonzero for all real numbers x.
\end{longlist}
\end{lemma}
\begin{pf}
If $\Pi$ is the L\'evy measure associated with the fragmentation
equation, then the assumption $\kappa<\infty$ is equivalent to $\int
_0^{\infty} x \Pi( d x)<\infty$, which, by Propositions 3.1
and 2.1 of \cite{CarmonaPetitYor} implies (i) and (ii). Next, it was
proven in the proof of Theorem 2 of \cite{HaasRegularity} that
$\mathbb E[I^{ix-1}] \neq0$ for all $x \in\mathbb R$ under the
additional assumption that the support of $\Pi$ is not included in a
set of the form $\{rn, n\geq0\}$ for some $r>0$.
\end{pf}
\begin{corollary}
\label{corofasterslower}
Suppose that (\ref{hypothese}) holds, $\kappa=\int_0^1 |{\ln(x)}| x
B( dx)<\infty$ and the support of $B$ is not included in a set
of the form $\{ a^n, n \in\mathbb N\}$ for some $a \in\,]0,1[$. Then,
for all measurable functions $g\dvtx[0,\infty[\  \rightarrow\,]0,\infty[ $
converging to $\rightarrow0$ at $\infty$,
\[
g(t)^{\alpha}\mathbb P \biggl( \biggl( \frac{\varphi(|\alpha
|t)}{|\alpha|t} \biggr)^{1/|\alpha|}\frac{X(t)}{g(t)} \leq1
\Big| X(t)>0 \biggr) \rightarrow\frac{ 1}{|\alpha|\kappa}.
\]
\end{corollary}
\begin{pf}
For $x \geq0$, let
\[
U_t(x):=g(t)^{\alpha}\mathbb P \biggl( \biggl(\frac{\varphi(|\alpha
|t)}{|\alpha|t} \biggr)^{1/|\alpha|}X(t) \leq xg(t) \Big|
X(t)>0 \biggr)
\]
and note that this quantity increases in $x$ when $t$ is fixed. Then
consider some random variable $I$, independent of $X$, with
distribution that of $\int_0^{\infty} \exp(\alpha\xi_r)\, d
r$. Consider $b$ such that $\mathbb P(I \leq b)>0$. Then
\[
U_t(x) \mathbb P(I \leq b) \leq g(t)^{\alpha}\mathbb P \biggl( \biggl(
\frac{\varphi(|\alpha|t)}{|\alpha|t} \biggr)^{1/|\alpha|}X(t)
I^{1/|\alpha|} \leq b^{1/|\alpha|} x g(t) \Big| X(t)>0
\biggr),
\]
which, according to Lemma \ref{Lemma2}(i), is bounded from above by
some constant (independent of $t$ and $x$) times $bx^{|\alpha|}$ for
all $x \geq0$ and $t \geq t_0$. That is, there exists some finite
constant $C$ such that for all $t$ sufficiently large and all $x \geq0$,
%
%
\begin{equation}
\label{majoU}
x^{\alpha}U_t(x) \leq C.
\end{equation}
Now, consider an increasing function $l\dvtx\mathbb N \rightarrow\mathbb
N$. For all $x \geq0$, the sequence $(U_{l(n)}(x), n \geq0)$ is
bounded. Hence, there exists some nondecreasing right-continuous
function $U \dvtx[0,\infty[\  \rightarrow[0,\infty[$, with $U(0)=0$, and
a subsequence $(U_{\tilde{l}(n)}, n \geq0)$ of $(U_{l(n)}, n \geq0)$
such that $U_{\tilde{l}(n)}(x) \rightarrow U(x) $ for a.e. $x>0$; see,
for example, \cite{Feller}, Theorem 2, Section VIII.7. Hence, if we
prove that the limit $U$ is given by
%
%
\begin{equation}
\label{caracU}
U(x)=\frac{x^{|\alpha|}}{ |\alpha| \kappa}\qquad \forall x\geq0,
\end{equation}
for all sequences $(l(n),n \geq0), (\tilde l(n),n \geq0)$ as defined
above, then we will have the expected result [note that the continuity
of the function involved in (\ref{caracU}) implies that the
convergence will hold for every $x \geq0$].

To prove (\ref{caracU}), recall that by Lemma \ref{lemmak}(ii), $\int
_0^{\infty} x^{-1}k(x) \, d x<\infty$. Hence, by dominated
convergence, for all $a>0$, $\int_0^{\infty} U_{l(n)}(ax^{1/\alpha
})k(x)\, d x \rightarrow\break\int_0^{\infty} U(ax^{1/\alpha})\times k(x)
\mathbb d x$. By Lemma \ref{Lemma2}(ii), we therefore have
%
%
\begin{equation}
\label{carac}
\int_0^{\infty} U(ax^{1/\alpha}) k(x) \, d x=a^{|\alpha|}\qquad
\forall a>0.
\end{equation}
We claim that this equation characterizes $U$ under the additional
assumption that the support of $B$ is not included in a set of the form
$\{ a^n, n \in\mathbb N\}$ for some $a \in\,]0,1[$. Indeed, note first
that by setting $V (x ):=\exp(x)U(\exp(x/\alpha))$ and
$\overline{k}{(x)}:=k(\exp(-x))$ for all $x \in\mathbb R$, the above
equation can be rewritten
\[
\int_{-\infty}^{\infty} V (x ) \overline{k} (y-x
) \, d x=1\qquad \forall y \in\mathbb R.
\]
However, the function $V$ is bounded a.e. on $\mathbb R$, by (\ref
{majoU}). Moreover, by Lem\-ma~\ref{lemmak}(ii), $\overline{k} \in
L^1(\mathbb R)$ and by Lemma \ref{lemmak}(iii), the Fourier transform
of $\overline{k}$ is nonzero on $\mathbb R$. We conclude, using the
Wiener approximation theorem for $L^1(\mathbb R)$ (\cite{BGT}, Theorem
4.8.4), that the above equation in $V$ has a unique bounded
solution (in the sense that two solutions are equal a.e.). This
determines $V$, hence $U$, almost everywhere. Since $U$ is
right continuous, it is determined for all $x \geq0$. Finally, it is
not hard to check that the expression for $U$ given by (\ref{caracU})
indeed satisfies (\ref{carac}).
\end{pf}

\subsection{Proof of Proposition
\protect\ref{fasterslower}\textup{(ii)}}
We need only prove the second part of Proposition \ref
{fasterslower}(ii), the first part being obvious since $\mu
_t(]1,\infty[)=0$ for all $t \geq0$ and $g(t) (\varphi(|\alpha
|t)/|\alpha|t)^{1/\alpha} \rightarrow\infty$ as $t \rightarrow
\infty$.
We keep the notation from the previous section and recall that we work
under the assumption (\ref{hypothese}). From the proof of Lemma \ref
{Lemma2}, we get that
\[
\ln\biggl( \frac{m ( t (1+ {|\alpha|h(t)^{|\alpha
|}}/{\varphi(|\alpha|t)} ) )}{m(t)} \biggr) \mathop{\sim
}_{t \rightarrow\infty} - h(t)^{|\alpha|}
\]
for all positive functions $h$ such that $h(t)^{|\alpha|}/\varphi(t)
\rightarrow0$ as $t \rightarrow\infty$. In other words, for such
functions $h$,
\[
\ln\biggl(m(t)^{-1}\mathbb P \biggl( \biggl( \frac{\varphi(|\alpha
|t)}{|\alpha|t} \biggr)^{1/|\alpha|}X(t) I^{1/|\alpha|} \geq
h(t) \biggr) \biggr) \mathop{\sim}_{t \rightarrow\infty} -
h(t)^{|\alpha|}.
\]
Note that for all $t \geq0$ and all $c>0$, since $X$ is independent of $I$,
\begin{eqnarray*}
&&\ln\biggl(m(t)^{-1}\mathbb P \biggl( \biggl( \frac{\varphi
(|\alpha|t)}{|\alpha|t} \biggr)^{1/|\alpha|}X(t)\geq\bigl(c\phi
\bigl(h(t)^{|\alpha|}\bigr) \bigr)^{1/|\alpha|} \biggr) \biggr)\\
&&\quad{} + \ln
\bigl( \mathbb P \bigl(I^{1/|\alpha|} \geq h(t)/ \bigl(c\phi\bigl(h(t)^{|\alpha
|}\bigr) \bigr)^{1/|\alpha|} \bigr) \bigr) \\
&&\qquad \leq\ln\biggl(m(t)^{-1}\mathbb P \biggl( \biggl( \frac{\varphi
(|\alpha|t)}{|\alpha|t} \biggr)^{1/|\alpha|}X(t) I^{1/|\alpha|}
\geq h(t) \biggr) \biggr).
\end{eqnarray*}
Suppose,\vspace*{1pt} moreover, that $h(t) \rightarrow\infty$ as $t \rightarrow
\infty$ and that $\beta<1$, which implies that $h(t)^{|\alpha|}/\phi
(h(t)^{|\alpha|}) \rightarrow\infty$. By Lemma \ref{Victor}, we
have, for all real numbers $c>0$,
\begin{eqnarray*}
-\ln\bigl( \mathbb P \bigl(I \geq h(t)^{|\alpha|}/c\phi
\bigl(h(t)^{|\alpha|}\bigr) \bigr) \bigr)
&\mathop{\sim}\limits_{\infty}&
\frac{1-\beta}{|\alpha|} \varphi\biggl(\frac{ |\alpha|h(t)^{|\alpha
|}}{c\phi(h(t)^{|\alpha|})} \biggr)\\
&\mathop{\sim}\limits_{\infty}&
\frac{(1-\beta)|\alpha|^{\beta/(1-\beta)}}{c^{1/(1-\beta)}}
h(t)^{|\alpha|},
\end{eqnarray*}
using both the regular variation of $\varphi$ and the fact that
$\varphi$ is the inverse of $t \rightarrow t/\phi(t)$ near $\infty$.
Now, let $\varepsilon\in\,]0,1[$ and $c$ be such that $c^{1/(1-\beta
)}>(1-\beta)|\alpha|^{\beta/(1-\beta)}$.
We have proven that
\begin{eqnarray*}
&&\ln\biggl(m(t)^{-1}\mathbb P \biggl( \biggl( \frac{\varphi(|\alpha
|t)}{|\alpha|t} \biggr)^{1/|\alpha|}X(t)
\geq
\bigl(c\phi \bigl(h(t)^{|\alpha|}\bigr) \bigr)^{1/|\alpha|} \biggr) \biggr)\\
&&\qquad\leq
-(1-\varepsilon) \biggl(1 - \frac{(1-\beta)|\alpha|^{\beta/(1-\beta
)}}{c^{1/(1-\beta)}} \biggr)h(t)^{|\alpha|}
\end{eqnarray*}
for $t$ large enough.
Next, let $g_{h,c}(t)= (c\phi(h(t)^{|\alpha|})
)^{1/|\alpha|}, t \geq0$, and suppose that $\beta>0$ (hence the
existence of the inverse of $\phi$ near $\infty$). We have, for $t$
large enough,
\begin{eqnarray*}
&&\ln\biggl(m(t)^{-1}\mathbb P \biggl( \biggl( \frac{\varphi(|\alpha|
t)}{|\alpha| t} \biggr)^{1/|\alpha|}X(t) \geq g_{h,c}(t) \biggr)
\biggr) \\
&&\qquad\leq -(1-\varepsilon) \biggl(1 - \frac{(1-\beta)|\alpha
|^{\beta/(1-\beta)}}{c^{1/(1-\beta)}} \biggr) \phi
^{-1}\bigl(g_{h,c}(t)^{|\alpha|}/c\bigr) \\
&&\qquad\mathop{\sim}_{\infty} -(1-\varepsilon) \biggl(1- \frac{(1-\beta
)|\alpha|^{\beta/(1-\beta)}}{c^{1/(1-\beta)}} \biggr) c^{-1/\beta
} \phi^{-1}\bigl(g_{h,c}(t)^{|\alpha|}\bigr).
\end{eqnarray*}
It is not hard to check that the maximum of
\[
\biggl\{ \biggl(1 - \frac{(1-\beta)|\alpha|^{\beta/(1-\beta
)}}{c^{1/(1-\beta)}} \biggr) c^{-1/\beta}, c>(1-\beta)^{1-\beta
}|\alpha|^{\beta} \biggr\}
\]
is equal to $\beta/ |\alpha|$ and is reached at $c=|\alpha|^{\beta}$.
Finally, if $g_h=g_{h, |\alpha|^{\beta}}$, letting $\varepsilon
\rightarrow0$, we have proven that
%
%
\begin{eqnarray}
\label{eq9}
&&\limsup_{t \rightarrow\infty} \frac{1}{\phi^{-1}(g_h(t)^{|\alpha
|})}\ln\biggl(m(t)^{-1}\mathbb P \biggl( \biggl( \frac{\varphi
(|\alpha| t)}{|\alpha| t} \biggr)^{1/|\alpha|}X(t) \geq g_h(t)
\biggr) \biggr) \nonumber\\[-8pt]\\[-8pt]
&&\qquad \leq-\frac{\beta} {|\alpha|}.\nonumber
\end{eqnarray}
To conclude, to get the second part of the statement of Proposition
\ref{fasterslower}(ii), suppose that $0<\beta<1$ and consider some
positive function $g$ that converges to $\infty$ at $\infty$, such
that $g(t)^{|\alpha|}t/\varphi(t) \rightarrow0$. Set $h(t)=(\phi
^{-1}(g(t)^{|\alpha|}/ |\alpha|^{\beta})^{1/|\alpha|}, t \geq0$.
Then, $h(t)$ converges to $\infty$ as $t \rightarrow\infty$ and it
is easily seen that $h(t)^{|\alpha|}/\varphi(t) \rightarrow0$ as $t
\rightarrow\infty$. Since $g=g_h$ with the notation above, the result
follows from (\ref{eq9}).

\section{Some properties of the limit measure
$ \mu_{\infty}$}
\label{PropertiesZ}

Recall that the distribution $x\mu_{\infty}( dx)$ on
$]0,\infty[$ is that of
$
R^{1/|\alpha|},
$
where $R$ denotes a random variable with entire positive moments
%
%
\begin{equation}
\label{Rmoments}
\mathbb E[R^n]=\phi(|\alpha|)\cdots\phi(n|\alpha|),\qquad n \geq1,
\end{equation}
that characterize its distribution. Using this particular moments'
shape, we get the following description of the measure $\mu_{\infty}$
near $0$ and $\infty$. Some of these properties are then used at the
end of this section to prove Proposition \ref{mass1}.
\begin{proposition}[(Behavior at $\infty$)]
\label{behaviorinfty}
\begin{longlist}
\item Suppose that (\ref{hypothese}) holds for some $\beta\in\,
]0,1[$. Then
\[
-\ln\biggl(\int_t^{\infty} x \mu_{\infty} ( d x) \biggr)
=- \ln\bigl( \mathbb P\bigl(R>t^{|\alpha|}\bigr) \bigr) \mathop{\sim
}_{\infty} \frac{\beta}{|\alpha|} \phi^{-1} \bigl(t^{|\alpha|}\bigr),
\]
where $\phi^{-1}$ denotes the inverse of $\phi$ (and is therefore a
function regularly varying at $\infty$ with index $1/\beta$).
\item Suppose that $\phi(\infty):=\int_0^{1}xB( d
x)<\infty$. Then $\mu_{\infty}$ has a bounded support with supremum
$\phi(\infty)^{1/|\alpha|}$
and
\[
\mu_{\infty} \bigl( \bigl\{\phi(\infty)^{1/|\alpha|} \bigr\} \bigr) >0
\quad\Leftrightarrow\quad\int^1 \frac{B( d x)}{1-x}<\infty.
\]
\end{longlist}
\end{proposition}
\begin{proposition}[(Behavior at $0$)]
Suppose that $\int_{0}^{\exp(-u)} x B( dx) $
varies regularly at $\infty$ with index $-\gamma, \gamma\in[0,1]$.
Then, as $s \rightarrow0$,
\[
\int_s^{\infty} x^{1+\alpha} \mu_{\infty} ( d x)=\mathbb
E\bigl[\mathbf1_{\{R>s^{|\alpha|} \}}R^{-1}\bigr]\sim\frac{1}{(|\alpha
|^{\gamma} \Gamma(1+\gamma)\phi(-1/\ln(s^{|\alpha|})))}
\]
and
\[
\int_0^s x \mu_{\infty}( d x)= \mathbb P\bigl(R<s^{|\alpha
|}\bigr)=\circ\biggl(\frac{s^{|\alpha|}}{\phi(-1/\ln(s^{|\alpha|}))}
\biggr).
\]
\end{proposition}
\begin{pf}
This is a direct consequence of Corollary 1 of Caballero and
Rivero \cite{CaballeroRivero}, which gives these results in terms of the
random variable $R$.
\end{pf}
\begin{pf*}{Proof of Proposition \protect\ref{behaviorinfty}} (i) Our
proof strongly relies on the proof of Proposition 2 of Rivero \cite
{VictorLog}. Rivero shows there that if a positive random variable $Y$
has entire moments satisfying
\[
\mathbb E[Y^n]=\prod_{i=1}^n\psi(i)
\]
for some function $\psi$ regularly varying at $\infty$ with index
$\gamma\in\,]0,1[$, then
\[
-\ln\bigl( \mathbb P(Y>t) \bigr)\mathop{\sim}_{\infty}\gamma\psi
^{\leftarrow}(t),
\]
where $\psi^{\leftarrow}$ is the right inverse of $\psi$. We apply
this result to the random variable $R$, by taking $\psi=\phi(|\alpha
| \cdot)$ and $\gamma=\beta$.

(ii) Using (\ref{Rmoments}) and the fact that $\phi$ is increasing,
we get, for all $n \geq0$,
%
%
\begin{equation}
\label{eq10}
\mathbb E \biggl[ \biggl(\frac{R}{\phi(\infty)} \biggr)^n \biggr]
\leq1.
\end{equation}
Besides, writing
\begin{eqnarray*}
\mathbb E \biggl[ \biggl(\frac{R}{\phi(\infty)} \biggr)^n \biggr]
&=& \mathbb E \biggl[ \biggl(\frac{R}{\phi(\infty)} \biggr)^n \mathbf
1_{\{ R>\phi(\infty)\}} \biggr] +\mathbb E \biggl[ \biggl(\frac
{R}{\phi(\infty)} \biggr)^n \mathbf1_{\{ R<\phi(\infty)\}}
\biggr]\\
&&{} + \mathbb P \bigl(R=\phi(\infty) \bigr)
\end{eqnarray*}
and using the monotone and dominated convergence theorems, we see that
\[
\mathbb E \biggl[ \biggl(\frac{R}{\phi(\infty)} \biggr)^n \biggr]
\mathop{\rightarrow}_{n \rightarrow\infty} \cases{
\infty, &\quad if $\mathbb P\bigl(R>\phi(\infty)\bigr)>0$,\cr
\mathbb P\bigl(R=\phi(\infty)\bigr), &\quad otherwise.}
\]
In particular, from (\ref{eq10}), we see that $\mathbb P(R>\phi
(\infty))=0$. Similarly, it is easy to show, using (\ref{Rmoments}),
that for all $0<\varepsilon<\phi(\infty)$,
\[
\mathbb E \biggl[ \biggl(\frac{R}{\phi(\infty)-\varepsilon}
\biggr)^n \biggr] \rightarrow\infty,
\]
which implies that $\mathbb P(R>\phi(\infty)-\varepsilon)>0$.
Finally, to get the remaining part of the statement, note that
\[
\ln\biggl( \frac{\phi(n|\alpha|)}{\phi(\infty)} \biggr) \mathop
{\sim}_{n \rightarrow\infty} \frac{\phi(n|\alpha|)}{\phi(\infty
)}-1=\frac{-1}{\phi(\infty)}\int_0^1 x^{n|\alpha|+1}B( dx).
\]
Therefore,
\[
\ln\biggl( \mathbb E \biggl[ \biggl(\frac{R}{\phi(\infty)}
\biggr)^n \biggr] \biggr)=\sum_{i=1}^{n} \ln\biggl( \frac{\phi(i|\alpha
|)}{\phi(\infty)} \biggr)
\]
converges to $-\infty$ as $n \rightarrow\infty$ if and only if
\[
\int_0^1 \sum_{i=1}^{\infty} x^{i|\alpha|+1}B( dx)=\int_0^1
\frac{x^{1+|\alpha|}}{1-x^{|\alpha|}}B( d x)=\infty.
\]
Since $\int_0xB( d x)<\infty$,
\[
\mathbb E \biggl[ \biggl(\frac{R}{\phi(\infty)} \biggr)^n \biggr]
\mathop{\rightarrow}_{n \rightarrow\infty}0 \quad\mbox{iff}\quad \int^1
\frac{1}{1-x}B( d x)=\infty,
\]
which ends the proof.
\end{pf*}
\begin{pf*}{Proof of Proposition \protect\ref{mass1}} From the
construction (\ref{solution}) of $\mu_t$, we see that
\[
\mu_t(\{ 1\})=\mu_0(\{ 1\}) \mathbb P \bigl(\xi(\rho(t))=0
\bigr)=\mu_0(\{ 1\})\mathbb P \bigl(\xi(t)=0 \bigr)
\]
and from the Poisson point process construction of a pure jump
subordinator with L\'evy measure $\Pi$, we have that $\mathbb P
(\xi(t)=0 )=\exp(-t\Pi(]0,\infty[))= \exp(-t\phi
(\infty) )$. Next, we get, from the factorization (\ref{face}), that
\[
\exp(-t\phi(\infty) )=\mathbb P \bigl( RI>t\phi(\infty
) \bigr) \geq\mathbb P(I>t) \mathbb P\bigl(R \geq\phi(\infty)\bigr).
\]
On the one hand, from the proof of Proposition \ref{behaviorinfty}, we
see that when $\phi(\infty)<\infty$, $\mathbb P(R \geq\phi(\infty
))=\mathbb P(R = \phi(\infty))$ and that this quantity is nonzero iff
$\int^1(1-x)^{-1}B( d x)<\infty$. On the other hand, under
(\ref{hypothese}), we get from the regular variation of $-\ln(\mathbb
P(I>t))$ that $\mathbb P(I>x^{\alpha}t)/\mathbb P(I>t) \rightarrow0$,
for all $0< x<1$, as $t \rightarrow\infty$ and then, from the
dominated convergence theorem that
\[
\frac{m(t)}{\mathbb P(I>t)}=\int_0^1 \frac{\mathbb P(I>x^{\alpha
}t)}{\mathbb P(I>t)} x \mu_0( d x) \rightarrow\mu_0(\{ 1\})
\qquad\mbox{as } t \rightarrow\infty.
\]
In other words, we have proven that under the hypothesis (\ref
{hypothese}), when $\mu_0(\{ 1\})>0$ and $\phi(\infty)<\infty$,
\[
\liminf_{t \rightarrow\infty} \frac{\mu_t(\{1\})}{m(t)} \geq
\mathbb P\bigl(R=\phi(\infty)\bigr)=\phi(\infty)^{1/|\alpha|} \mu_{\infty
}\bigl(\bigl\{\phi(\infty)^{1/|\alpha|}\bigr\}\bigr).
\]
Next, suppose that (\ref{hypothese}) holds, $\int_0 |{\ln(x)}| x
B( dx)<\infty$ and $\phi(\infty)<\infty$. According to
Theorem \ref{theq}, for all $\varepsilon\in\,]0,1[$ such that
$(1-\varepsilon)\phi(\infty)^{1/|\alpha|}$ is not an atom of $\mu
_{\infty}$,
\[
\frac{\mu_t(\{1\})}{m(t)} \leq\frac{\int_{1-\varepsilon}^1x\mu
_t( dx)}{m(t)} \mathop{\rightarrow}_{t \rightarrow\infty}
\int_{(1-\varepsilon) \phi(\infty)^{1/|\alpha|}}^{\phi(\infty
)^{1/|\alpha|}}x\mu_{\infty}( d x).
\]
Letting $\varepsilon\rightarrow0$, we get
\[
\limsup_{t \rightarrow\infty} \frac{\mu_t(\{1\})}{m(t)}\leq\phi
(\infty)^{1/|\alpha|} \mu_{\infty}\bigl(\bigl\{\phi(\infty)^{1/|\alpha|}\bigr\}\bigr).
\]
\upqed\end{pf*}

\section{Examples}
\label{Examples}

Below is a list of standard examples where the main quantities involved
in our results can be computed explicitly. More precisely, for each of
these examples, we specify the distributions of $I$ [defined in (\ref
{defI})] and $R$ [defined in (\ref{defmuR})], which leads to explicit
expressions of the limit measure $\mu_{\infty}$ [since $R^{1/|\alpha
|} \stackrel{\dd}\sim x\mu( dx)$] and of the mass
\[
m_1(t)=\mathbb P(I>t),\qquad t \geq0,
\]
which is the mass of the solution of the fragmentation equation
starting from \mbox{$\mu_0=\delta_1$}.
We also specify the behavior as $t \rightarrow\infty$ of the quantity
$\varphi(|\alpha| t)/|\alpha| t$ involved in the statement of
Theorem \ref{theq}. For all of these examples, we give the main tools
to get the distributions of $I$ and $R$, but we leave the calculation
details to the reader. We recall that $\beta$ denotes the index of
regular variation of the hypothesis (\ref{hypothese})
and that when $(\mu_t,t \geq0)$ is a solution of the equation with
parameters $(\alpha,B)$, $(\mu_{ct},t \geq0)$ is a solution of the
equation with parameters $(\alpha,cB)$. For this reason, in the
examples below, given a measure $B$, we choose its ``representative''
among the measures $cB$, $c>0$, which is the most convenient for the
statement of the results.

The first four examples concern absolutely continuous measures
$B( du)=b(u) \, du$, where $b$ is a function defined on
$]0,1[$. The L\'evy measure is therefore also absolutely continuous and
we denote by $\pi$ its density. It turns out that the limit
distribution $\mu_{\infty}$ is also absolutely continuous. We denote
by $u_{\infty}$ its density.
\begin{Example}[{[$b(u)=bu^{b-2}$, $b>0$;
$\alpha<0$]}]\label{Example1}
\begin{itemize}
\item $\beta=0;$
\item $\varphi(t) \sim t$ as $t \rightarrow\infty;$
\item $I\stackrel{\dd}{\sim} \Gamma(b/{|\alpha|}+1,1);$
\item $m_1(t)=\frac{1}{\Gamma(b/|\alpha|+1)} \int
_t^{\infty} x^{b/|\alpha|}\exp(-x) \, dx$, $t \geq0;$
\item $R \stackrel{\dd}{\sim} \beta
(1,b/{|\alpha|});$
\item $u_{\infty}(x)=b x^{|\alpha|-2} ( 1-x^{|\alpha
|} )^{{b}/{|\alpha|}-1}$, $0<x<1$.
\end{itemize}
The notation $\Gamma(x,y)$ [resp., $\beta(x,y)]$ refers to the
classical Gamma distribution with parameters $x,y>0$ (resp., Beta
distribution). In these examples, the density of the L\'evy measure
associated with $B$ is $
\pi(x)=b \exp(-bx)$, $x >0$, hence the L\'evy measure associated with
the subordinator $|\alpha|\xi$ (where $\xi$ has L\'evy measure $\Pi
)$ has a density given by $b \exp(-bx/|\alpha|)/|\alpha|$, $x >0$.
According to Example B, page 5 of \cite{CarmonaPetitYor}, the
density of $I$ is then proportional to $x^{b/|\alpha|} \exp(-x),
x>0$. Finally, we refer to formula (4), Section 3 of \cite{BYFacExp},
to get the distribution of $R$.

We point out that the solutions of the fragmentation equation with this
measure $B$ are studied in \cite{McGZ}. In particular, when $\alpha
=-b/2$ and $\mu_0=\delta_1$, the solutions $(\mu_t, t \geq0)$ have
the explicit expression
\[
\mu_t( dx)=\exp(-t) \bigl( \delta_1( d x)+ b x^{b-2}
\bigl(t-\tfrac{1}{2}t^2(1-x^{-b/2}) \bigr) \mathbf
1_{\{0<x<1\}}\, d x \bigr),
\]
which gives
\[
m_1(t)=\exp(-t) \biggl(1+t+\frac{t^2}{2} \biggr)
\]
and, for all bounded test functions $f\dvtx]0,\infty[\  \rightarrow
\mathbb R$,
\[
\frac{1}{m_1(t)}\int_0^1 f(x) x \mu_t( dx) \mathop
{\rightarrow}_{t \rightarrow\infty} \int_0^1
f(x)bx^{b-1}(x^{-b/2}-1) \, d x,
\]
which is consistent with the above expressions for $m_1$ and $u_{\infty}$.
\end{Example}
\begin{Example}[{[$b(u)=|\alpha|\Gamma(1-\gamma)^{-1} u^{{|\alpha|}/{\gamma}-2}(1-u^{{|\alpha
|}/{\gamma}} )^{-\gamma-1}$; $0<\gamma<1$; \mbox{$\alpha<0$}]}]\label{Example2}

\begin{itemize}
\item $\beta=\gamma$;
\item $\varphi(t) \sim(\frac{\gamma}{|\alpha
|} )^{{\gamma}/({1-\gamma})} t^{{1}/({1-\gamma})}$ as $t
\rightarrow\infty$;
\item $I\stackrel{\dd}{\sim} \tau_{\gamma
}^{-\gamma} $;
\item $m_1(t)=\int_t^{\infty} g_{\gamma}(x) \,
dx$, $t >0$;
\item $ R \stackrel{\dd}{\sim} \mathbf
e(1)^{\gamma}$;
\item $\mu_{\infty}(x)=|\alpha|\gamma^{-1}x^{
{|\alpha|}/{\gamma}-2}\exp(-x^{|\alpha|/\gamma} )$, $x>0$.
\end{itemize}
Here, $\mathbf e(1)$ denotes a random variable with exponential
distribution with parameter 1 and $\tau_{\gamma}$ denotes a $\gamma
$-stable random variable, that is, with Laplace transform $t \in
[0,\infty[\  \rightarrow\exp(-t^\gamma)$. Hence, $\tau_{\gamma
}^{-\gamma} $ has the so-called Mittag--Leffler distribution. We
recall that it possesses a density given by
\[
g_{\gamma}(x)=\frac{1}{\pi\gamma} \sum_{i=0}^{\infty} \frac
{(-x)^{i-1}}{i !} \Gamma(\gamma i +1) \sin(\pi\gamma i),\qquad x >0,
\]
and its entire positive moments are equal to $n!/\Gamma(\gamma n +1)$,
$\forall n \geq1$ (see, e.g., \cite{PitmanStFl}, Section 0.3).
The L\'evy measure associated with $B$ has a density for $x>0$ given by
\[
\pi(x)=\frac{|\alpha|\exp(-|\alpha|x/\gamma)}{\Gamma(1-\gamma
)(1-\exp(-|\alpha|x/\gamma))^{\gamma+1}}.
\]
Using formula (5) and the following discussion in \cite{BYFacExp}, we
get that $I\stackrel{\dd}{\sim} \tau_{\gamma}^{-\gamma} $
and $ R \stackrel{\dd}{\sim} \mathbf e(1)^{\gamma}$.
\end{Example}
\begin{Example}[{[$b(u)=|\alpha|\gamma^2 ((1- \gamma)\Gamma(2-\gamma) )^{-1} u^{{\gamma|\alpha|}/({1-\gamma})
-2}(1-\break u^{|\alpha|/(1-\gamma)})^{-\gamma-1}$;
$0<\gamma<1; \alpha<0$]}]\label{Example3}

\begin{itemize}
\item $\beta=\gamma$;
\item $\varphi(t) \sim(1-\gamma)^{-1}|\alpha|^{
{\gamma}/({\gamma-1})}t^{{1}/({1-\gamma})}$, as $t \rightarrow
\infty$;
\item $I \stackrel{\dd}{\sim} \mathbf
e(1)^{1-\gamma}$;
\item $m_1(t) = \exp(-t^{1/(1-\gamma)})$, $t \geq0$;
\item $R \stackrel{\dd}{\sim} \tau_{1-\gamma
}^{\gamma-1}$;
\item $\mu_{\infty}(x)=|\alpha| x^{|\alpha
|-2}g_{1-\gamma} (x^{|\alpha|})$, $x>0$,
\end{itemize}
where $g_{1-\gamma}$ is the Mittag--Leffler density given in the
previous example. Note the duality with this previous example. In the
present example,
\[
\pi(x)=\frac{|\alpha|\gamma^2\exp(|\alpha|x/(1-\gamma
))}{(1-\gamma)\Gamma(2-\gamma)(\exp(|\alpha|x/(1-\gamma
))-1)^{1+\gamma}},\qquad x>0,
\]
and we again refer to formula (5) and the discussion which follows in
\cite{BYFacExp} to get the distributions of $I$ and $R$.
\end{Example}
\begin{Example}[{[$b(u)=|\alpha| \Gamma
(2+\alpha)^{-1} u^{|\alpha|-2}(1-u)^{\alpha-1}$; $-1 < \alpha<0$]}]\label{Example4}

\begin{itemize}
\item $\beta=|\alpha|$;
\item $\varphi(t) \sim( \frac{t}{1+\alpha}
)^{{1}/({1-|\alpha|})}$ as $t \rightarrow\infty$;
\item $I /(1+\alpha)$ is a size-biased version of the
Mittag--Leffler distribution with parameter $|\alpha|$, that is, for
all test functions $f$,
\[
\mathbb E[f(I)]=\frac{\mathbb E [f ((1+\alpha)\tau
_{|\alpha|}^{-|\alpha|} ) \tau_{|\alpha|}^{-|\alpha|}
]}{\mathbb E [\tau_{|\alpha|}^{-|\alpha|} ]};
\]
\item $m_1(t)=\Gamma(|\alpha|+1)\int_{t/(1+\alpha
)}^{\infty}xg_{|\alpha|}(x) \, dx$, $t >0$;
\item $ ((1+\alpha)R )^{1/|\alpha
|}\stackrel{\dd}{\sim} \Gamma(|\alpha|,1)$;
\item $\mu_{\infty}(x)=\frac{(1+\alpha)}{\Gamma
(|\alpha|)} x^{|\alpha|-2}\exp(-(1+\alpha)^{1/|\alpha|}x)$, $x>0$.
\end{itemize}
Indeed, here,
\[
\pi(x)=\frac{|\alpha| \exp(x)}{\Gamma(2+\alpha) (\exp(x)
-1)^{1-\alpha}},\qquad x>0.
\]
Following the end of the proof of Lemma 4 of Miermont \cite
{Miermont1}, we get that $I$ has its moment of order $k$ equal to
\[
\frac{k! (1+\alpha)^k\Gamma(|\alpha|)}{\Gamma( (k+1) |\alpha
| ) }
\]
for all $k \in\mathbb N$. Hence,
\begin{eqnarray*}
\mathbb E[R^k] &=& \frac{k!}{\mathbb E[I^k]}\\
&=&\frac{\Gamma( (k+1)
|\alpha| )} {(1+\alpha)^k\Gamma(|\alpha|)}\\
&=& \frac
{1}{(1+\alpha)^k\Gamma(|\alpha|)} \int_0^{\infty} x^{(k+1)|\alpha
|-1}\exp(-x) \, d x.
\end{eqnarray*}
\end{Example}
\begin{Remark*}
Note that Examples \ref{Example2}, \ref{Example3} and \ref{Example4} give, for all
$0<\gamma<1$:

\begin{itemize}
\item if $b(u)=u^{-1}(1-u)^{-\gamma-1}$ and $\alpha=-\gamma$,
then
\[
x\mu_{\infty}( dx)\stackrel{\dd} \sim c_1(\gamma)
\mathbf e(1);
\]

\item if $b(u)=u^{\gamma-2}(1-u)^{-\gamma-1}$ and $\alpha=\gamma-1$,
then
\[
x\mu_{\infty}( dx) \stackrel{\dd} \sim c_2(\gamma)\tau
_{1-\gamma}^{-1};
\]

\item if $b(u)=u^{\gamma-2}(1-u)^{-\gamma-1}$ and $\alpha=-\gamma$,
then
\[
x\mu_{\infty}( dx) \stackrel{\dd} \sim c_3(\gamma
)\Gamma(\gamma,1),
\]
\end{itemize}
where $c_1(\gamma), c_2(\gamma)$ and $c_3(\gamma)$ are real numbers
that depend on $\gamma$.
Hence, both $\alpha$ and the behavior of $b$ near $0$ have a
significant influence on the shape of the limit measure $\mu_{\infty}$.

Finally, we turn to the case where $B$ is a Dirac measure.
\end{Remark*}
\begin{Example}[{($B=a^{-1}\delta_{a}$ for some $a \in\,]0,1[$; $\alpha<0$)}]\label{Example5}
\begin{itemize}
\item $\beta=0$;
\item $\varphi(t) \sim t$ as $t \rightarrow\infty$;
\item $I$ has a density $k$ on $]0,\infty[$ given by
\[
k(x)=\sum_{i \geq0} \exp\bigl(\alpha\ln(a) i- x \exp(\alpha\ln
(a) i) \bigr) \prod_{p \neq i} \bigl(1-\exp\bigl(\alpha\ln(a)(i-p)\bigr)
\bigr)^{-1};
\]
\item $m_1(t)=\int_t^{\infty} k(x) \, d x$, $t \geq0$.
\end{itemize}
In this case,
$\Pi=\delta_{-\ln(a)}$, that is, the associated subordinator is a
Poisson process. We then refer to \cite{CarmonaPetitYor}, Proposition
6.5(ii), for the expression of the density $k$.
Note that $\phi(t)=(1-a^t)$ for all $t \geq0$, hence $\mathbb
E[R^n]=\prod_{i=1}^n (1-(|\alpha|a)^i)$ for all $n \geq1$.\vspace*{-10pt}
\end{Example}

\begin{appendix}\label{App}
\section*{Appendix: Existence and uniqueness of solutions}
This appendix is devoted to the proof of Theorem \ref{theq} on the
existence and uniqueness of solutions of the fragmentation equation
(\ref{eqfrag}). Therefore, in this section, $\alpha\in\mathbb R$.
The proof follows the main lines of that of Theorem 1 in \cite
{HaasLossMass}, which gives existence and uniqueness of solutions of a
slightly restricted form of the fragmentation equation (\ref{eqfrag})
and which concentrates on solutions starting from $\mu_0=\delta_1$.
We note that it was implicit in the statement of this theorem that a
solution should satisfy assumptions (\ref{gainmass}) and (\ref{creationmass}).

Let $\xi$ denote a subordinator with L\'evy measure $\Pi$ and zero
drift, such that \mbox{$\xi_0=0$}. We recall that its semigroup possesses the
Feller property and that the domain of its infinitesimal generator
contains at least all functions $f$ that are continuously
differentiable on $\mathbb R$ and such that $f$ and $f'$ tend to $0$ at
infinity; see, for example, Chapter 1 of \cite{BertoinLevy}. As a
consequence, the domain of the infinitesimal generator of $\exp(-\xi
)$ contains continuously differentiable functions $f$ on $[0,\infty[$
with compact support and null near 0.

One can easily check that when $f \dvtx[0,\infty[\  \rightarrow\mathbb R$
is bounded and continuous, the function
\[
x \rightarrow\mathbb E \bigl[f\bigl(x\exp\bigl(-\xi_{\rho(x^{\alpha
}t)}\bigr)\bigr) \bigr]
\]
is also bounded and continuous on $[0, \infty[$. This mainly relies on
the c\`adl\`ag and quasi-left-continuity (\cite{BertoinLevy},
Proposition 7, Chapter 1) of subordinators.

Now, for every $0 <a< b$, let $\mathcaligr C_{a,b}$ be the set of
continuous functions $f \dvtx[0,b] \rightarrow\mathbb R$ that are null on
$[0,a]$, and let $\mathcaligr C^1_{a,b}$ be the set of continuously
differentiable functions $f \dvtx[0,b] \rightarrow\mathbb R$ that are
null on $[0,a]$. It is clear from the remark above that for all
$0<a<b$, the linear operators $T_t$ and $\tilde{T}_t$, $t \geq0$,
defined by
\[
T_t(f)(x)=\mathbb E [f(x\exp(-\xi_t)) ]
\]
and
\[
\tilde{T}_t(f)(x)=\mathbb E \bigl[f\bigl(x\exp\bigl(-\xi_{\rho(x^{\alpha
}t)}\bigr)\bigr) \bigr]
\]
send $\mathcaligr C_{a,b}$ into $\mathcaligr C_{a,b}$.
Following the proof of Theorem 1 of \cite{HaasLossMass} (see also
\cite{Lamperticaract}), we see that both families of operators define
strongly continuous contraction semigroups on $\mathcaligr C_{a,b}$ and
that the domains of their infinitesimal generators are identical and
contain $\mathcaligr C^1_{a,b}$. These generators are, respectively, given,
for $f \in\mathcaligr C^1_{a,b}$ and $x \in[0,b]$, by
\[
A(f)(x)=\int_0^\infty\bigl( f(x\exp(-y))-f(x) \bigr) \Pi(
d y)
\]
and
\[
\tilde A(f)(x)=x^{\alpha}A(f)(x),\qquad x>0, \tilde A(f)(0)=0.
\]
Note that when $B$ is a measure on $]0,1[$ defined from $\Pi$ by (\ref
{PiB}), we have
\[
\tilde A(f)(x)=x^{\alpha} \int_0^1 \bigl(f(xy)-f(x) \bigr)y
B( dy).
\]

\subsection*{Existence of solutions to (\protect\ref{eqfrag})} With the
above remarks and Kolmogorov's backward equation (see Proposition 15,
page 9 of \cite{EK}), we have that
%
%
\begin{equation}
\label{Kolmogorov}
\tilde{T}_t(f)(x)=f(x)+\int_0^t \tilde{T}_s(\tilde{A}(f))(x)
\, ds,
\end{equation}
$\forall x \in[0,b]$, $\forall f \in\mathcaligr C^1_{a,b}$, $\forall
0<a<b$, $\forall b>0$. In other words, if we let $f\dvtx[0,\infty
[\ \rightarrow\mathbb R$ be null near $0$ and continuously
differentiable, then, considering its restriction to $[0,b]$ and $x
\leq b$, we have that $f$ and $x$ satisfy (\ref{Kolmogorov}).

Now, consider $\nu_0$, a probability measure on $]0,\infty[$, and
define for all $t>0$ a measure $\nu_t$ on $]0,\infty[$ by
\[
\langle\nu_t,g\rangle:=\langle\nu_0, \tilde{T}_t(g)\rangle
\]
for all bounded, measurable functions $g$ on $[0,\infty[$ such that
$g(0)=0$. Note that for all $t \geq0$, $\nu_t(]0,\infty[) \leq1$
and $\nu_t(x\geq M)=0$ provided that $\nu_0(x\geq M)=0$ for some $M>0$.
Then let $f$ be some continuously differentiable function on $[0,\infty
[$, null near 0 and with compact support. It is clear that $\tilde A
(f)$ is null near 0 and it is easy to see, using Fubini's theorem, that
there exist some constants $b,c>0$ such that $|\tilde A (f)| (x)\leq c
x^{\alpha}\overline{\Pi}(\ln(x/b))$ for large enough $x$ [here,
$\overline\Pi(y)=\int_y^{\infty}\Pi( dx)$]. In particular,
$\tilde A (f)$ is bounded on $[0,\infty[$ when $x^{\alpha} \overline
\Pi(\ln(x))$ is bounded near $\infty$ (hence when $\alpha\leq0$).
It is then clear that in such a case, we can apply Fubini's theorem
when integrating (\ref{Kolmogorov}) with respect to $\nu_0$ to get
\[
\langle\nu_t,f\rangle=\langle\nu_0,f\rangle+\int_0^t \langle\nu_s,
\tilde{A}(f)\rangle\, ds.
\]
This holds for all continuously differentiable functions $f$ on
$[0,\infty[$, null near 0 and with compact support. Therefore,
defining the measures $\mu_t$ on $]0,\infty[$ by $\langle\mu_t,g\rangle
:=\langle\nu
_t,\tilde g\rangle$, where $g$ denotes any test function on $]0,\infty
[$ and
$\tilde g(x)=g(x)/x$, $x>0$, we have proven that $(\mu_t,t\geq0)$ is
a solution of the fragmentation equation, as defined in the
\hyperref[sec1]{Introduction}. To summarise: provided that the function $x \rightarrow
x^{\alpha} \overline\Pi(\ln(x))$ is bounded near $\infty$, for all
measures $\mu_0$ on $]0,\infty[$ such that $\int_0^{\infty} x \mu
_0( d x)=1$, there exists a solution, constructed via
subordinators, of the fragmentation equation.\looseness=1

When $\alpha>0$, the function $x \rightarrow x^{\alpha} \overline\Pi
(\ln(x))$ may not be bounded near $\infty$. Another way to tackle the
problem in this case is to use the definition of $\rho$ to get that
\begin{eqnarray*}
&&\int_0^{\infty} \int_0^t \tilde{T}_s(|\tilde{A}(f)|)(x) \, ds
\nu_0( d x)\\
&&\qquad =\int_0^{\infty} \mathbb E \biggl[ \int_0^{\rho
(x^{\alpha}t)} |A(f)|(x\exp(-\xi_u)) \, d u \biggr] \nu
_0( d x),
\end{eqnarray*}
the function $f$ still being supposed continuously differentiable on
$[0,\infty[$, null near 0 and with compact support. For such $f$, the
function $A(f)$ is bounded on $[0,\infty[$. Hence, the double integral
involved in the identity above is bounded by a constant times $\int
_0^{\infty} \mathbb E [\rho(x^{\alpha}t) ] \nu_0(
d x)$, which is finite provided that $\int^{\infty} \ln(x) \nu
_0( dx)<\infty$: indeed, according to Proposition 2 of \cite
{BYFacExp}, for all $x\geq0$, $\mathbb E[\rho(x)]=\int_0^x \mathbb
E[\exp(-s\times R)] \, ds$, where $R$ is a random variable with
distribution $\mu_R$ defined by (\ref{defmuR}). Now, let $I$ be
a random variable defined by (\ref{defI}), independent of $R$, and
consider a real number $a$ such that $\mathbb P(I \leq a)>0$. Using the
factorization property (\ref{face}), we get
\begin{eqnarray*}
\mathbb E[\rho(x)] \mathbb P(I \leq a) &=& \int_0^x \mathbb E\bigl[\exp(-sR)
\mathbf1_{ \{I \leq a \}}\bigr] \, ds\\
 &\leq&\int_0^x \mathbb E[\exp(-sa^{-1} \mathbf e(1))] \, ds= a\ln(1+a^{-1}x).
\end{eqnarray*}
It this then possible to apply Fubini's theorem when integrating (\ref
{Kolmogorov}) with respect to $\nu_0$ and we conclude, as above, that
there exists a solution of (\ref{eqfrag}).

\subsection*{Uniqueness of solutions of (\protect\ref{eqfrag})}
Let $\nu_0$ be a probability measure with support included in $]0,b]$
for some $b>0$ and suppose that $(\nu_t,t \geq0)$ is a family of
measures with support included in $]0,b]$ such that
\[
\langle\nu_t,f\rangle=\langle\nu_0,f\rangle+\int_0^t \langle\nu
_s,\tilde A(f)\rangle\, d s,\qquad
\forall f \in\bigcup_{0<a<b} \mathcaligr C^1_{a,b}.
\]
Suppose, moreover, that $\nu_t(]0,\infty[) \leq1 $, $\forall t \geq
0$. Our goal is to prove that $(\nu_t,t \geq0)$ is uniquely
determined. Using the fact that the total weight of $\nu_t$ is less
than or equal to 1, we get that $\sup_{t \geq0} \langle\nu_t, |f
|\rangle<\infty
$ and $\sup_{t \geq0} \langle\nu_t, |A(f) |\rangle<\infty$ for each $f
\in\bigcup
_{0<a<b} \mathcaligr C^1_{a,b}$.
It is then possible to follow the proof of Proposition 18, Section 4.9
of \cite{EK} to deduce that uniqueness holds, provided that, for all
$\lambda>0$, $(\lambda\,\mathrm{id} -\tilde A(f))(\mathcaligr C^1_{a,b})$
is dense (for the uniform norm) in $\mathcaligr C_{a,b}$ for all $0<a<b$.
Following the proof of Theorem 1 in \cite{HaasLossMass}, we see that
$\mathcaligr C^1_{a,b}$ is a core for the strongly contraction semigroup $
\tilde{T_t} \dvtx\mathcaligr C_{a,b} \rightarrow\mathcaligr C_{a,b} $, $t
\geq
0$. Hence the result.
\end{appendix}

\section*{Acknowledgments}
The author wishes to warmly thank
St\'ephane Mischler, who inspired this paper by asking whether a
probabilistic approach could be used to study the large time behavior
of solutions to the fragmentation equation with shattering. Many thanks
also to Victor Rivero for pointing out the reference \cite{KM} to
generalize his Proposition 2 in \cite{VictorLog}.

%

%
\printaddresses


\begin{thebibliography}{30}

\bibitem{Aldous96}
%
\begin{bincollection}[mr]
\bauthor{\bsnm{Aldous},~\bfnm{David}\binits{D.}}
(\byear{1996}).
\btitle{Probability distributions on cladograms}.
In \bbooktitle{Random Discrete Structures ({M}inneapolis, {MN}, 1993)}.
\bseries{The IMA Volumes in Mathematics and its Applications}
\bvolume{76}
\bpages{1--18}.
\bpublisher{Springer}, \baddress{New York}.
\bid{mr={1395604}}
\end{bincollection}
%
\endbibitem

\bibitem{Banasiaksurvey}
%
\begin{barticle}[mr]
\bauthor{\bsnm{Banasiak},~\bfnm{J.}\binits{J.}}
(\byear{2006}).
\btitle{Shattering and non-uniqueness in fragmentation models---an analytic
approach}.
\bjournal{Phys. D}
\bvolume{222}
\bpages{63--72}.
\bid{mr={2265768}}
\end{barticle}
%
\endbibitem

\bibitem{BertoinLevy}
%
\begin{bbook}[mr]
\bauthor{\bsnm{Bertoin},~\bfnm{Jean}\binits{J.}}
(\byear{1996}).
\btitle{L\'evy Processes}.
\bseries{Cambridge Tracts in Mathematics}
\bvolume{121}.
\bpublisher{Cambridge Univ. Press}, \baddress{Cambridge}.
\bid{mr={1406564}}
\end{bbook}
%
\endbibitem

\bibitem{BertoinAB}
%
\begin{barticle}[mr]
\bauthor{\bsnm{Bertoin},~\bfnm{Jean}\binits{J.}}
(\byear{2003}).
\btitle{The asymptotic behavior of fragmentation processes}.
\bjournal{J. Eur. Math. Soc. (JEMS)}
\bvolume{5}
\bpages{395--416}.
\bid{mr={2017852}}
\end{barticle}
%
\endbibitem

\bibitem{BertoinBook}
%
\begin{bbook}[mr]
\bauthor{\bsnm{Bertoin},~\bfnm{Jean}\binits{J.}}
(\byear{2006}).
\btitle{Random Fragmentation and Coagulation Processes}.
\bseries{Cambridge Studies in Advanced Mathematics}
\bvolume{102}.
\bpublisher{Cambridge Univ. Press}, \baddress{Cambridge}.
\bid{mr={2253162}}
\end{bbook}
%
\endbibitem

\bibitem{BC}
%
\begin{barticle}[mr]
\bauthor{\bsnm{Bertoin},~\bfnm{Jean}\binits{J.}} \AND
\bauthor{\bsnm{Caballero},~\bfnm{Maria-Emilia}\binits{M.-E.}}
(\byear{2002}).
\btitle{Entrance from {$0+$} for increasing semi-stable {M}arkov processes}.
\bjournal{Bernoulli}
\bvolume{8}
\bpages{195--205}.
\bid{mr={1895890}}
\end{barticle}
%
\endbibitem

\bibitem{BMEnergy}
%
\begin{barticle}[mr]
\bauthor{\bsnm{Bertoin},~\bfnm{Jean}\binits{J.}} \AND
\bauthor{\bsnm{Mart{\'{\i}}nez},~\bfnm{Servet}\binits{S.}}
(\byear{2005}).
\btitle{Fragmentation energy}.
\bjournal{Adv. in Appl. Probab.}
\bvolume{37}
\bpages{553--570}.
\bid{mr={2144567}}
\end{barticle}
%
\endbibitem

\bibitem{BYFacExp}
%
\begin{barticle}[mr]
\bauthor{\bsnm{Bertoin},~\bfnm{Jean}\binits{J.}} \AND
\bauthor{\bsnm{Yor},~\bfnm{Marc}\binits{M.}}
(\byear{2001}).
\btitle{On subordinators, self-similar {M}arkov processes and some
factorizations of the exponential variable}.
\bjournal{Electron. Comm. Probab.}
\bvolume{6}
\bpages{95--106}.
\bid{mr={1871698}}
\end{barticle}
%
\endbibitem

\bibitem{BGT}
%
\begin{bbook}[mr]
\bauthor{\bsnm{Bingham},~\bfnm{N.~H.}\binits{N.~H.}},
\bauthor{\bsnm{Goldie},~\bfnm{C.~M.}\binits{C.~M.}} \AND
\bauthor{\bsnm{Teugels},~\bfnm{J.~L.}\binits{J.~L.}}
(\byear{1989}).
\btitle{Regular Variation}.
\bseries{Encyclopedia of Mathematics and Its Applications}
\bvolume{27}.
\bpublisher{Cambridge Univ. Press}, \baddress{Cambridge}.
\bid{mr={1015093}}
\end{bbook}
%
\endbibitem

\bibitem{CaballeroRivero}
%
\begin{barticle}[mr]
\bauthor{\bsnm{Caballero},~\bfnm{Maria~Emilia}\binits{M.~E.}} \AND
\bauthor{\bsnm{Rivero},~\bfnm{V{\'{\i}}ctor~M.}\binits{V.~M.}}
(\byear{2009}).
\btitle{On the asymptotic behaviour of increasing self-similar {M}arkov
processes}.
\bjournal{Electron. J. Probab.}
\bvolume{14}
\bpages{865--894}.
\bid{mr={2497455}}
\end{barticle}
%
\endbibitem

\bibitem{CarmonaPetitYor}
%
\begin{bincollection}[mr]
\bauthor{\bsnm{Carmona},~\bfnm{Philippe}\binits{P.}},
\bauthor{\bsnm{Petit},~\bfnm{Fr{\'e}d{\'e}rique}\binits{F.}} \AND
\bauthor{\bsnm{Yor},~\bfnm{Marc}\binits{M.}}
(\byear{1997}).
\btitle{On the distribution and asymptotic results for exponential functionals
of {L}\'evy processes}.
In \bbooktitle{Exponential Functionals and Principal Values Related to
{B}rownian Motion}
\bpages{73--130}.
\bpublisher{Rev. Mat. Iberoam.}, \baddress{Madrid}.
\bid{mr={1648657}}
\end{bincollection}
%
\endbibitem

\bibitem{ChowCuzick}
%
\begin{barticle}[mr]
\bauthor{\bsnm{Chow},~\bfnm{Y.~S.}\binits{Y.~S.}} \AND
\bauthor{\bsnm{Cuzick},~\bfnm{Jack}\binits{J.}}
(\byear{1979}).
\btitle{Moment conditions for the existence and nonexistence of optimal
stopping rules for {$S\sb{n}/n\sp{1}$}}.
\bjournal{Proc. Amer. Math. Soc.}
\bvolume{75}
\bpages{300--307}.
\bid{mr={532155}}
\end{barticle}
%
\endbibitem

\bibitem{ESRR}
%
\begin{barticle}[mr]
\bauthor{\bsnm{Escobedo},~\bfnm{M.}\binits{M.}},
\bauthor{\bsnm{Mischler},~\bfnm{S.}\binits{S.}} \AND
\bauthor{\bsnm{Rodriguez~Ricard},~\bfnm{M.}\binits{M.}}
(\byear{2005}).
\btitle{On self-similarity and stationary problem for fragmentation and
coagulation models}.
\bjournal{Ann. Inst. H. Poincar\'e Anal. Non Lin\'eaire}
\bvolume{22}
\bpages{99--125}.
\bid{mr={2114413}}
\end{barticle}
%
\endbibitem

\bibitem{EK}
%
\begin{bbook}[mr]
\bauthor{\bsnm{Ethier},~\bfnm{Stewart~N.}\binits{S.~N.}} \AND
\bauthor{\bsnm{Kurtz},~\bfnm{Thomas~G.}\binits{T.~G.}}
(\byear{1986}).
\btitle{Markov Processes: Characterization and Convergence}.
\bpublisher{Wiley}, \baddress{New York}.
\bid{mr={838085}}
\end{bbook}
%
\endbibitem

\bibitem{Feller}
%
\begin{bbook}[mr]
\bauthor{\bsnm{Feller},~\bfnm{William}\binits{W.}}
(\byear{1971}).
\btitle{An Introduction to Probability Theory and Its Applications}
\bvolume{2},
\bedition{2nd} ed.
\bpublisher{Wiley}, \baddress{New York}.
\bid{mr={0270403}}
\end{bbook}
%
\endbibitem

\bibitem{Filippov}
%
\begin{barticle}[auto:springertagbib-v1.0]
\bauthor{\bsnm{Filippov},~\bfnm{A.}\binits{A.}}
(\byear{1961}).
\btitle{On the distribution of the sizes of particles which undergo splitting}.
\bjournal{Theory Probab. Appl.}
\bvolume{6}
\bpages{275--294}.
\end{barticle}
%
\endbibitem

\bibitem{FG}
%
\begin{barticle}[mr]
\bauthor{\bsnm{Fournier},~\bfnm{Nicolas}\binits{N.}} \AND
\bauthor{\bsnm{Giet},~\bfnm{Jean-S{\'e}bastien}\binits{J.-S.}}
(\byear{2003}).
\btitle{On small particles in coagulation-fragmentation equations}.
\bjournal{J. Stat. Phys.}
\bvolume{111}
\bpages{1299--1329}.
\bid{mr={1975930}}
\end{barticle}
%
\endbibitem

\bibitem{FGdensities}
%
\begin{barticle}[mr]
\bauthor{\bsnm{Fournier},~\bfnm{Nicolas}\binits{N.}} \AND
\bauthor{\bsnm{Giet},~\bfnm{Jean-S{\'e}bastien}\binits{J.-S.}}
(\byear{2006}).
\btitle{Existence of densities for jumping stochastic differential equations}.
\bjournal{Stochastic Process. Appl.}
\bvolume{116}
\bpages{643--661}.
\bid{mr={2205119}}
\end{barticle}
%
\endbibitem

\bibitem{HaasLossMass}
%
\begin{barticle}[mr]
\bauthor{\bsnm{Haas},~\bfnm{B{\'e}n{\'e}dicte}\binits{B.}}
(\byear{2003}).
\btitle{Loss of mass in deterministic and random fragmentations}.
\bjournal{Stochastic Process. Appl.}
\bvolume{106}
\bpages{245--277}.
\bid{mr={1989629}}
\end{barticle}
%
\endbibitem

\bibitem{HaasRegularity}
%
\begin{barticle}[mr]
\bauthor{\bsnm{Haas},~\bfnm{B{\'e}n{\'e}dicte}\binits{B.}}
(\byear{2004}).
\btitle{Regularity of formation of dust in self-similar fragmentations}.
\bjournal{Ann. Inst. H. Poincar\'e Probab. Statist.}
\bvolume{40}
\bpages{411--438}.
\bid{mr={2070333}}
\end{barticle}
%
\endbibitem

\bibitem{KM}
%
\begin{barticle}[mr]
\bauthor{\bsnm{K{\"o}nig},~\bfnm{Wolfgang}\binits{W.}} \AND
\bauthor{\bsnm{M{\"o}rters},~\bfnm{Peter}\binits{P.}}
(\byear{2002}).
\btitle{Brownian intersection local times: Upper tail asymptotics and thick
points}.
\bjournal{Ann. Probab.}
\bvolume{30}
\bpages{1605--1656}.
\bid{mr={1944002}}
\end{barticle}
%
\endbibitem

\bibitem{Lamperticaract}
%
\begin{barticle}[mr]
\bauthor{\bsnm{Lamperti},~\bfnm{John}\binits{J.}}
(\byear{1972}).
\btitle{Semi-stable {M}arkov processes. {I}}.
\bjournal{Z. Wahrsch. Verw. Gebiete}
\bvolume{22}
\bpages{205--225}.
\bid{mr={0307358}}
\end{barticle}
%
\endbibitem

\bibitem{LaurencotCoagFrag}
%
\begin{barticle}[mr]
\bauthor{\bsnm{Lauren{\c{c}}ot},~\bfnm{Philippe}\binits{P.}}
(\byear{2000}).
\btitle{On a class of continuous coagulation-fragmentation equations}.
\bjournal{J. Differential Equations}
\bvolume{167}
\bpages{245--274}.
\bid{mr={1793195}}
\end{barticle}
%
\endbibitem

\bibitem{McGZ}
%
\begin{barticle}[mr]
\bauthor{\bsnm{McGrady},~\bfnm{E.~D.}\binits{E.~D.}} \AND
\bauthor{\bsnm{Ziff},~\bfnm{Robert~M.}\binits{R.~M.}}
(\byear{1987}).
\btitle{``{S}hattering'' transition in fragmentation}.
\bjournal{Phys. Rev. Lett.}
\bvolume{58}
\bpages{892--895}.
\bid{mr={927489}}
\end{barticle}
%
\endbibitem

\bibitem{Melzak}
%
\begin{barticle}[mr]
\bauthor{\bsnm{Melzak},~\bfnm{Z.~A.}\binits{Z.~A.}}
(\byear{1957}).
\btitle{A scalar transport equation}.
\bjournal{Trans. Amer. Math. Soc.}
\bvolume{85}
\bpages{547--560}.
\bid{mr={0087880}}
\end{barticle}
%
\endbibitem

\bibitem{Miermont1}
%
\begin{barticle}[mr]
\bauthor{\bsnm{Miermont},~\bfnm{Gr{\'e}gory}\binits{G.}}
(\byear{2003}).
\btitle{Self-similar fragmentations derived from the stable tree. {I}.
{S}plitting at heights}.
\bjournal{Probab. Theory Related Fields}
\bvolume{127}
\bpages{423--454}.
\bid{mr={2018924}}
\end{barticle}
%
\endbibitem

\bibitem{PitmanStFl}
%
\begin{bbook}[mr]
\bauthor{\bsnm{Pitman},~\bfnm{J.}\binits{J.}}
(\byear{2006}).
\btitle{Combinatorial Stochastic Processes}.
\bseries{Lecture Notes in Math.}
\bvolume{1875}.
\bpublisher{Springer}, \baddress{Berlin}.
\bid{mr={2245368}}
\end{bbook}
%
\endbibitem

\bibitem{VictorLog}
%
\begin{barticle}[mr]
\bauthor{\bsnm{Rivero},~\bfnm{V{\'{\i}}ctor}\binits{V.}}
(\byear{2003}).
\btitle{A law of iterated logarithm for increasing self-similar {M}arkov
processes}.
\bjournal{Stoch. Stoch. Rep.}
\bvolume{75}
\bpages{443--472}.
\bid{mr={2029617}}
\end{barticle}
%
\endbibitem

\bibitem{Wagner}
%
\begin{barticle}[mr]
\bauthor{\bsnm{Wagner},~\bfnm{Wolfgang}\binits{W.}}
(\byear{2005}).
\btitle{Explosion phenomena in stochastic coagulation-fragmentation models}.
\bjournal{Ann. Appl. Probab.}
\bvolume{15}
\bpages{2081--2112}.
\bid{mr={2152254}}
\end{barticle}
%
\endbibitem

\bibitem{ZMcG}
%
\begin{barticle}[mr]
\bauthor{\bsnm{Ziff},~\bfnm{Robert~M.}\binits{R.~M.}} \AND
\bauthor{\bsnm{McGrady},~\bfnm{E.~D.}\binits{E.~D.}}
(\byear{1985}).
\btitle{The kinetics of cluster fragmentation and depolymerisation}.
\bjournal{J. Phys. A}
\bvolume{18}
\bpages{3027--3037}.
\bid{mr={814641}}
\end{barticle}
%
\endbibitem

\end{thebibliography}
\end{document}